\documentclass[11pt]{article}
\usepackage{amsmath,amssymb,theorem}
\usepackage[dvips]{graphicx}

\newtheorem{theorem}{Theorem}[section]
\newtheorem{proposition}[theorem]{Proposition}
\newtheorem{corollary}[theorem]{Corollary}
\newtheorem{lemma}[theorem]{Lemma}
\newtheorem{remark}[theorem]{Remark}
\newtheorem{definition}[theorem]{Definition}
\newtheorem{pr}{Proof.}

\newcommand{\proof}[1]{
\begin{pr}
#1 \,\fbox{}
\end{pr}
}

\begin{document}

\title{Lam\'{e} operators with projective octahedral and icosahedral monodromies}
\author{Keiri Nakanishi }
\date{November 18, 2004}
\maketitle

\begin{abstract}

We show that there exists a Lam\'{e} operator $L_{n}$ with projective octahedral monodromy for each 
$n \in \tfrac{1}{2}( \mathbb{N} + \tfrac{1}{2}) \cup \tfrac{1}{3}( \mathbb{N} + \tfrac{1}{2}) $, 
and with projective icosahedral monodromy for each 
$n \in \tfrac{1}{3}( \mathbb{N} + \tfrac{1}{2}) \cup \tfrac{1}{5}( \mathbb{N} + \tfrac{1}{2}) $.
To this end, we construct Grothendieck's dessins d'enfants corresponding 
to the Belyi morphisms which pull-back hypergeometric operators into Lam\'{e} operators 
$L_{n}$ with the desired monodromies.

\end{abstract}

\section{Introduction}
In this paper, we consider Lam\'{e} differential operators
\[ L_{n} = \left(\frac{d}{d x} \right) ^{2} + \frac{1}{2} \left(\frac{1}{x} + \frac{1}{x - 1} + 
\frac{1}{x - \lambda}\right)\cdot \frac{d}{d x} - \frac{n(n + 1)x + B}{4x(x - 1)(x - \lambda)} \]
whose solutions are algebraic over $\mathbb{C}(x)$.
Here, $\lambda$ is a complex number with $\lambda \neq 0, 1$, 
$n$ (so-called degree parameter) is a rational number, 
and $B \in \mathbb{C}$ is the accessory parameter.
The possible finite projective monodromies of $L_{n}$ were studied by Baldassarri, Chiarellotto, 
and Dwork, and recently by Beukers, Dahmen, Li\c{t}canu, van der Waall, and Zapponi.
One of the most remarkable results is that there are at most finitely many equivalence classes 
of $L_{n}$ for fixed $n \in \mathbb{Q}$ and fixed finite monodromy group.
This was first done by Chiarellotto~\cite{C}, and later shown by Li\c{t}canu by using the 
notion of Grothendieck's dessins d'enfants~\cite{L1}. 
For details of the theory of Grothendieck's dessins d'enfants, see~\cite{S} and~\cite{SV}. 
Moreover, Chiarellotto and Li\c{t}canu got the explicit formula for the number of equivalence 
classes with projective dihedral monodromy of order $2N$ for the case $n = 1$, 
which has been generalized more recently by Dahmen~\cite{D} for arbitrary $n$. 
They translated the counting problem for the number of equivalence classes of Lam\'{e} operators 
into that for the number of the dessins compatible with the ramification data of Belyi morphisms 
which pull-back hypergeometric operators into the Lam\'{e} operators. 
This strategy is based on the Klein's theorem, 
which claims that a second order Fuchsian differential operator with finite projective monodromy 
is a rational pull-back of hypergeometric operator in the ``\textit{basic Schwarz list}". 

To carry out this program, the method of Grothendieck's dessins d'enfants 
by Litcanu and Dahmen provides a powerful tool. 
Baldassarri determined the possible finite projective monodromy groups of $L_{n}$~\cite{B2}, 
but recently, Li\c{t}canu got the same results and the necessary conditions for $n$ to have 
fixed possible finite projective monodromy group by using the notion of Grothendieck's dessins 
d'enfants~\cite{L2}. 
By~\cite{B2} and~\cite{L2}, the possible finite projective monodormy groups are dihedral group $\mathcal{D}_{2N}$, 
octahedral group $\mathcal{S}_{4}$, or icosahedral group $\mathcal{A}_{5}$. 
More recently, Li\c{t}canu~\cite{L2} proves the following theorem:

\vspace{1ex}
\noindent
\begin{theorem}[~\cite{B2}, ~\cite{L2} Theorem 3.4]\label{Litcanu}

\hspace{2mm}
$(1)$ If the projective monodromy group of Lam\'{e} operator $L_{n}$ is dihedral,
then $n \in \mathbb{Z}$

\noindent
$(2)$ If the projective monodromy group of Lam\'{e} operator $L_{n}$ is octahedral,
then $n \in \tfrac{1}{2}( \mathbb{Z} + \tfrac{1}{2}) \cup \tfrac{1}{3}( \mathbb{Z} + 
\tfrac{1}{2}) $.

\noindent
$(3)$ If the projective monodromy group of Lam\'{e} operator $L_{n}$ is icosahedral,
then $n \in \tfrac{1}{3}( \mathbb{Z} + \tfrac{1}{2}) \cup \tfrac{1}{5}( \mathbb{Z} + 
\tfrac{1}{2}) $.

\noindent
$(4)$ There is no Lam\'{e} operator with projective cyclic monodromy.

\noindent
$(5)$ There is no Lam\'{e} operator with projective tetrahedral monodromy.
\end{theorem}

\vspace{1ex}
The proof of this theorem is based on the analysis of the Belyi morphism which pull-backs 
the hypergeometric operator into the Lam\'{e} operator, as well as the combinatorial 
data of the corresponding dessin. 
Conversely, the following problem arises: 

\vspace{1ex}
\noindent
\underline{\textbf{Problem}}

\noindent
(1) \textit{For each} $n \in \mathbb{Z}$, 
\textit{does there exist Lam\'{e} operator} $L_{n}$ 
\textit{with projective dihedral monodromy?}

\noindent
(2) \textit{For each} $n \in \tfrac{1}{2}(\mathbb{Z} + \tfrac{1}{2}) \cup 
\tfrac{1}{3}(\mathbb{Z} + \tfrac{1}{2})$, 
\textit{does there exist Lam\'{e} operator} $L_{n}$ 
\textit{with projective octahedral monodromy?}

\noindent
(3) \textit{For each} $n \in \tfrac{1}{3}(\mathbb{Z} + \tfrac{1}{2}) \cup 
\tfrac{1}{5}(\mathbb{Z} + \tfrac{1}{2})$, 
\textit{does there exist Lam\'{e} operator} $L_{n}$ 
\textit{with projective icosahedral monodromy?}

\vspace{1ex}
If we replace $n$ by $- n - 1$, it is easy to see that $L_{n} = L_{- n - 1}$. 
Hence we can assume $n > - \tfrac{1}{2}$. 
As we saw, (1) is solved by Beukers and van der Waall~\cite[Theorem5.1]{BW} and Dahmen~\cite{D}. 
Beukers and van der Waall~\cite[Theorem 6.1]{BW} and Baldassarri~\cite[(3.e)]{B2} 
gave some examples of Lam\'{e} operators which have projective 
octahedral monodormy and icosahedral monodromy. 

\vspace{1ex}
The aim of this article is to solve (2) and (3) of the problem above. 
Assuming $n>-\tfrac{1}{2}$, we have a few possible negative $n$ in each case.
Such ``exceptional'' cases will be dealt with case by case in Remarks \ref{rem-exc1} and \ref{rem-exc2}, where we will see that these cases can be easily dismissed. 
Thus we may assume $n\geq 0$.
In this situation, the following theorem gives the existence of Lam\'e operators $L_n$ with projective octahedral monodromy and projective icosahedral monodromy for each $n$ as in Theorem \ref{Litcanu} (2) and (3):

\vspace{1ex} 
\underline{\textbf{Main Theorem}}

\noindent
(1) \textit{For each} $n \in \tfrac{1}{2}(\mathbb{N} + \tfrac{1}{2}) 
\cup \tfrac{1}{3}(\mathbb{N} + \tfrac{1}{2})$, 
\textit{there exists a Lam\'{e} operator} $L_{n}$ \textit{with projective octahedral monodromy.}

\noindent
(2) \textit{For each} $n \in \tfrac{1}{3}(\mathbb{N} + \tfrac{1}{2}) 
\cup \tfrac{1}{5}(\mathbb{N} + \tfrac{1}{2})$, 
\textit{there exists a Lam\'{e} operator} $L_{n}$ \textit{with projective icosahedral monodromy.}

\vspace{1ex}
We will prove this theorem by constructing explicitly the dessins compatible with the ramification data of 
the Belyi morphisms which pull-back hypergeometric operators with the same projective 
monodromy group into the Lam\'{e} operators $L_{n}$ for each $n$.  
From this theorem, we can see that there exist infinitely Lam\'{e} operators with projective 
octahedral monodromy and infinitely many ones with projective icosahedral monodromy, which seems unknown. 
Note that this theorem does not answer the counting problem of the numbers of the equivalence 
classes of Lam\'{e} operators with projective octahedral and icosahedral monodromies. 

\section{Preliminaries.}

Our first aim in this section is to reduce the existences of the Lam\'{e} operators with 
projective octahedral (resp. icosahedral) monodromy to the existences of the 
Belyi morphisms which pull-back the hypergeometric operators with projective octahedral 
(resp. icosahedral) monodromies into the Lam\'{e} operators. 
The second aim is to reduce them to the existences of the corresponding dessins. 

\subsection{Hypergeometric operators and Lam\'{e} operators.}

In this subsection, we review some results on hypergeometric operators, their rational 
pull-backs, and Lam\'{e} operators. 

\vspace{1ex}
Let us first consider the linear differential operator on $\mathbb{P}^{1}$: 
\begin{equation}
L = D^{n} + a_{1}(z) \cdot D^{n-1} + \cdots + a_{n-1}(z) \cdot D + a_{n}(z), 
\label{eq:L}
\end{equation}
with $D^{i} = (\tfrac{d}{dz})^{i}$ and 
$a_{j}(z)$ rational functions in $\mathbb{C}(z)$ 
for $1 \leq i, j \leq n$. 
The linear operator (\ref{eq:L}) is said to be $Fuchsian$ if 
any point on $\mathbb{P}^{1}$ is regular or regular singular. 

Throughout this paper, we treat projective monodromies of the Fuchsian operators 
rather than (full) monodromies of them. 
Let us consider the natural projection $P$ : 
\[
P : GL(n, \mathbb{C}) \longrightarrow PGL(n, \mathbb{C}). 
\]
The monodromy group $G$ of the operator (\ref{eq:L}) is defined in $GL(n, \mathbb{C})$. 
Then its natural image of $G$ by $P$ is said to be the \textit{projective monodromy} of operator (\ref{eq:L}), 
and we denote it by $PG$, 
i.e., $PG = G \cdot Z / Z$ 
where $Z = \{\lambda \cdot I_{n} | \lambda \in \mathbb{C}^{\ast} \}$. 
Here, $PG$ is a subgroup of $PGL(n, \mathbb{C})$ and determined up to conjugate. 

\vspace{1ex}
Let us consider a second order Fuchsian operator with finite projective monodromy. 
If it has precisely three regular singular points, it is 
the so-called hypergeomtric operator and it has the following normalized form
\[
H_{ \lambda, \mu, \nu} = \left(\frac{d}{dx} \right)^{2} + \left\{ \frac{1 - \lambda ^{2}}{4x^{2}} 
+ \frac{1 - \mu ^{2}}{4(x - 1)^{2}} + \frac{\lambda ^{2} + \mu ^{2} + \nu ^{2} - 1}{4x(x - 1)} \right \}, 
\]
where $\lambda + \mu + \nu > 1$. 
The regular singular points of $H_{\lambda, \mu, \nu}$ are $0, 1, \infty$, and 
their exponent differences are $\lambda, \mu, \nu,$ respectively. 
The finite projective monodromy groups of $H_{\lambda, \mu, \nu}$ 
are classified as in the following ``\textit{basic Schwarz list}". 

\begin{center}
  \begin{tabular}{|c|c|}                                                                                \hline
    \multicolumn{1}{|c|}{$(\lambda, \mu, \nu )$}      & projective monodromy of $H_{\lambda, \mu, \nu}$ \\ \hline 
     $(1/n, 1,1/n)$                &         $\mathcal{C}_{n}$ : cyclic of order $n$ \\ \hline
     $(1/2, 1/n, 1/2)$     &     $\mathcal{D}_{2n}$ : dihedral of order $2n$ \\ \hline
     $(1/2, 1/3, 1/3)$     &                 $\mathcal{A}_{4}$ : tetrahedral \\ \hline
     $(1/2, 1/3, 1/4)$     &                  $\mathcal{S}_{4}$ : octahedral \\ \hline
     $(1/2, 1/3, 1/5)$     &                 $\mathcal{A}_{5}$ : icosahedral \\ \hline
  \end{tabular}
\end{center}

In general, second order Fuchsian operators with finite projective monodromy are characterized by 
the following theorem by Klein. 

\begin{theorem}[Klein]
Let $L$ be a second order Fuchsian operator with finite projective monodromy $PG$ in normalized form on 
$\mathbb{P}^{1}$, i.e. $L = (\tfrac{d}{dx})^{2} + Q(x)$, $Q(x) \in \mathbb{C}(x)$.
Then there exists a morphism $f : \mathbb{P}^{1} \rightarrow \mathbb{P}^{1}$ 
which ramifies at most over the set $\{0, 1,\infty \}$, 
and a unique hypergeometric operator $H$ in the Schwarz list, having the same projective monodromy $PG$, 
such that $f^{\ast} H = L$. Moreover, the morphism $f$ as above is unique up to M\"{obius} transformations 
except in the case 
$(\lambda, \mu, \nu) = (1/2, 1/2, 1/2)$.    
\end{theorem}

\proof{
\rm{
See~\cite{K} or~\cite[Theorem 1.8]{B1}.
}
}

Let $C$ be an algebraic curve defined over $\mathbb{C}$.
A morphism $f : C \rightarrow \mathbb{P}^{1}$ is said to be  \textit{Belyi morphism} 
if $f$ has at most three critical values. 

\vspace{1ex}
\textit{A Lam\'{e} operator} is a second order Fuchsian operator 
having four regular singular points on $\mathbb{P}^{1}$ 
with exponent differences $\tfrac{1}{2}, \tfrac{1}{2}, \tfrac{1}{2},$ and $|n + \tfrac{1}{2}|$ 
where $n$ is a rational number. 
If its four regular singular points are $0, 1, \lambda,$ and $\infty$, 
and their exponent differences are $\tfrac{1}{2}, \tfrac{1}{2}, \tfrac{1}{2},$ 
and $|n + \tfrac{1}{2}|$ respectively, 
then, after suitable transformation, we can assume it has the following Riemann scheme : 

\[
  \begin{Bmatrix} 
     x = 0 & x = 1 & x = \lambda & x = \infty \\ 
     0 & 0 & 0 & - \tfrac{n}{2} \\ 
     \tfrac{1}{2} & \tfrac{1}{2} & \tfrac{1}{2} & \tfrac{n+1}{2}  
  \end{Bmatrix}, 
\]

\vspace{1ex}
\noindent
and the Lam\'{e} operator $L_{n}$ has the following form: 

\[ 
L_{n} = \left(\frac{d}{d x} \right) ^{2} + \frac{1}{2} \left(\frac{1}{x} + \frac{1}{x - 1} + 
\frac{1}{x - \lambda}\right)\cdot \frac{d}{dx} - \frac{n(n + 1)x + B}{4x(x - 1)(x - \lambda)}, 
\]
where $B \in \mathbb{C}$ is the accessory parameter. 
Therefore, if there exists a Belyi morphism $f : \mathbb{P}^{1} \rightarrow \mathbb{P}^{1}$ satisfying 
the following condition ($\bigstar$) bellow, then $f^{\ast}H_{\lambda, \mu, \nu}$ is a Lam\'{e} operator. 

\vspace{1ex}
\textbf{Condition} ($\bigstar$) : 
$f^{\ast}H_{\lambda, \mu, \nu}$ has four regular singular points $0, 1, \lambda,$ and $\infty$ 
and their exponent differences are 
$\tfrac{1}{2}, \tfrac{1}{2}, \tfrac{1}{2},$ and $|n + \tfrac{1}{2}|$ respectively.

\vspace{1ex}

We summarize some facts about pull-backs of Fuchsian operators. 
We refer to~\cite[Ch.2]{vdW} for details. 

\begin{proposition} \label{prop1}
Let $L$ be a Fuchsian operator on $\mathbb{P}^{1}$, 
and $f : \mathbb{P}^{1} \rightarrow \mathbb{P}^{1}$ be a morphism. 
Then $f^{\ast} L$ is again Fuchsian.
\end{proposition}

\proof{
\rm{
See~\cite[Proposition 2.6.3]{vdW}.
}
}

\begin{proposition}\label{prop2}
Let $L$ be a Fuchsian operator on $\mathbb{P}^{1}$ with projective monodromy $PG_{L}$, 
and $f^{\ast}L$ is a pull-back by a morphism 
$f : \mathbb{P}^{1} \rightarrow \mathbb{P}^{1}$, with projective monodrory $PG_{f^{\ast}L}$ . 
Then $PG_{f^{\ast}L}$ is conjugate in $PGL(2, \mathbb{C})$ to a subgroup of $PG$. 
\end{proposition}

\proof{
\rm{
See~\cite[Corollary 2.6.10]{vdW}.
}
}

\begin{corollary}\label{cor3} 
Let $f : \mathbb{P}^{1} \rightarrow \mathbb{P}^{1}$ be a Belyi morphism satisfying the condition 
$($$\bigstar$$)$, 
$H_{1/2, 1/3, 1/4}$ $($resp. $H_{1/2, 1/3, 1/5}$$)$ 
the hypergeometric operator with projective octahedral $($resp. icosahedral$)$ monodromy. 
Let 
$n \in \tfrac{1}{2}(\mathbb{Z} + \tfrac{1}{2}) 
\cup \tfrac{1}{3}(\mathbb{Z} + \tfrac{1}{2}) 
\cup \tfrac{1}{5}(\mathbb{Z} + \tfrac{1}{2})$. 
Then the second order Fuchsian operator $f^{\ast}H_{1/2, 1/3, 1/4}$ 
$($resp. $f^{\ast}H_{1/2, 1/3, 1/5}$$)$ have projective octahedral 
$($resp. icosahedral$)$ monodromy. 
\end{corollary}

\proof{
\rm{
By Theorem \ref{Litcanu}, the projective monodromy of $f^{\ast}H_{1/2, 1/3, 1/4}$ 
is octahedral or icosahedral, and by Proposition \ref{prop2}, 
it must be conjugate to a subgroup of projective octahedral group. 
But projective icosahedral group cannot be conjugate to a subgroup of projective 
octahedral group, then the projective monodromy of $f^{\ast}H_{1/2, 1/3, 1/4}$ 
is octahedral. The icosahedral case is proved similarly. 
}
}

By this Corollary, the construction of Lam\'{e} operator with projective octahedral
monodromy amounts to the construction of the Belyi morphism satisfying the condition ($\bigstar$).

\subsection{Belyi morphisms and Grothendieck's dessins d'enfants.}

This subsection gives some reviews about \textit{Grothendieck's dessins d'enfants}. 
For more details, we refer to~\cite{S} and~\cite{SV}. 
Let us first recall Belyi's Theorem. 

\begin{theorem}[Belyi's Theorem]
Let $X$ be an algebraic curve over $\mathbb{C}$. 
Then $X$ is defined over $\bar{\mathbb{Q}}$ if and only if there exist a morphism 
$\beta : X \rightarrow \mathbb{P}^{1}(\mathbb{C})$ which ramifies at most over $\{0, 1, \infty \}$.  
\end{theorem}

\proof{
\rm{
Well-known; see~\cite{Be} or~\cite[Theorem I.2]{S}. 
}
}

Let $\beta : X \rightarrow \mathbb{P}^{1}$ be a Belyi morphism. 
For a point $P \in X$, we denote by $e_{P}$ the ramification index at $P$ of $\beta$. 
The Belyi morphism $\beta$ is said to be \textit{clean} 
if $e_{P} = 2$ for any $P \in \beta^{-1}(1)$, 
and \textit{preclean} if $e_{P} \leq 2$ for any $P \in \beta^{-1}(1)$. 
Consider a pair $(X, \beta)$ consisting of a complex algebraic curve 
defined over $\bar{\mathbb{Q}}$ 
and a morphism $\beta : X \rightarrow \mathbb{P}^{1}$. 
The pair $(X, \beta)$ is said to be a \textit{Belyi pair} 
if the morphism $\beta$ ramifies at most over $\{0, 1, \infty \}$. 
Two pairs $(X, \beta)$ and $(Y, \alpha)$ are said to be \textit{isomorphic} 
if there exists an isomorphism $\phi : X \rightarrow Y$ 
such that $\beta = \alpha \circ \phi$.

\begin{definition}[Dessins d'enfants]
\rm{
Let $X$ be a compact Riemann surface, 
$X_{1}$ a connected 1-complex, 
$X_{0}$ the set of vertices of $X_{1}$, 
$[\iota]$ an isotopical class of inclusions $\iota : X_{1} \hookrightarrow X$. 
The triple $D = ( X_{0} \subset X_{1}, [\iota] )$ is said to be 
\textit{Grothendieck's dessin d'enfant} on $X$ 
if $D$ satisfies the following conditions: 

\noindent 
(1) The complement of $X_{0}$ in $X_{1}$ is a finite disjoint union of segments 
and each segment is homeomorphic to the interval $(0, 1)$. 

\noindent
(2) The complement of $\iota (X_{1})$ in $X$ is a finite disjoint union of open cells
(simply connected regions). 

\noindent
(3) Each element of $X_{0}$ is equipped with the mark $``\bullet"$ or $``\ast"$ 
and if two different elements of $X_{0}$ are connected by a segment, 
one is equipped with $``\bullet"$ and another $``\ast"$. 

}
\end{definition}

\begin{definition}
\rm{
Two Grothendieck's dessins $D = (X_{0} \subset X_{1}, [\iota])$ on $X$ 
and $D' = (X'_{0} \subset X'_{1},[\iota '])$ on $X'$ 
are said to be \textit{equivalent} if there exists a homeomorphism 
$\phi : X \rightarrow X'$ such that 
$\phi |_{\iota(X_{1})} : \iota(X_{1}) \rightarrow \iota'(X'_{1})$ and 
$\phi |_{\iota(X_{0})} : \iota(X_{0}) \rightarrow \iota'(X'_{0})$ are homeomorphisms. 
}
\end{definition}

\begin{definition}
\rm{
A Grothendieck's dessin $D = (X_{0} \subset X_{1},[\iota] )$ 
is said to be \textit{preclean} if all vertices with the mark $``\ast"$ 
have valencies $\leq 2$. 
If all vertices with the mark $``\ast"$ have valencies 2, 
$D$ is said to be \textit{clean}. 
}
\end{definition}

\vspace{1ex}
Let $(X, \beta)$ be a Belyi pair. 
Then from the Belyi pair $(X, \beta)$, we can construct a dessin 
$D = (\beta^{-1}(\{0, 1 \}) \subset \beta^{-1}([0, 1]))$ 
by putting the mark $``\bullet"$ on the vertices of $\beta^{-1}(0)$, 
and $``\ast"$ on the vertices of $\beta^{-1}(1)$.

\begin{theorem}[Grothendieck Correspondence]
This correspondence gives a bijection between 
the set of isomorphic classes of preclean Belyi pairs
and the set of equivalence classes of preclean dessins. 
\end{theorem}

\proof{
\rm{
See~\cite[Theorem I.5]{S}. 
}
}

The procedure for getting Belyi pairs from dessins is given in~\cite[Chapter I, \S3]{S}. 
By this correspondence, the most important thing is that 
the ramification multiplicities of points in $\beta^{-1}(0)$ 
(resp. $\beta^{-1}(1)$) 
are translated to the valencies of $``\bullet"$ (resp. $``\ast"$). 

\vspace{1ex}
By Grothendieck Correspondence, we can construct Lam\'{e} operator $L_{n}$ with projective 
octahedral and icosahedral monodromy, if there exists a dessin d'enfant corresponding to a Belyi 
morphism which satisfies the condition ($\bigstar$). 
In the next section, we are going to construct the dessins that the corresponding Belyi 
morphisms satisfies the condition ($\bigstar$) for each $n$.

\begin{definition}

\rm{
A Belyi morphism $f : \mathbb{P}^{1} \rightarrow \mathbb{P}^{1}$ 
is said to be \textit{ $\ast$-morphism } 
if $\{0, 1,\infty \} \subseteq f^{-1} ( \{0, 1,\infty \})$.
}
\end{definition}

\begin{remark}
\rm{
Under the action of $PGL(2, \mathbb{C})$, any Belyi morphism 
$f : \mathbb{P}^{1} \rightarrow \mathbb{P}^{1}$ 
is transformed to a $\ast$-morphism. 
}
\end{remark}

\section{Constructions of the Dessins.}

\subsection{The case of projective octahedral monodromy.} 

We start this subsection by preparing some notations. 
We denote a second order Fuchsian differential operator on $\mathbb{P}^{1}$ by $L$ 
and its exponent difference at $P \in \mathbb{P}^{1}$ by $\Delta _{P, L}$. 
Set  
$\Delta _{L} = \displaystyle \sum_{P \in \mathbb{P}^{1}} ( \Delta_{P, L} - 1 ) $. 

\vspace{1ex} 
We need the following useful lemma. 

\begin{lemma}[~\cite{BD} Lemma 1.5]\label{Bal-Dw}
Let $f : \mathbb{P}^{1} \rightarrow \mathbb{P}^{1}$ be a morphism, 
and $L$ a Fuchsian second order differential operator. 
Then 
\[ 
\mathrm{deg}(f) = \frac{\Delta_{f^{\ast} L} + 2}{\Delta_{L} + 2}. 
\] 
\end{lemma}

\proof{
\rm{
This is the genus 0 case of~\cite[Lemma 1.5]{BD}. 
Let $Q \in \mathbb{P}^{1}$, $P \in \mathbb{P}^{1}$ be points with $f(Q) = P$. 
If $\alpha_{1}, \alpha_{2}$ are local exponents of $L$ at $P$, 
then the local exponents of $f^{\ast}L$ at $Q$ are $\alpha_{1} \cdot e_{Q, f}$ and 
$\alpha_{2} \cdot e_{Q, f}$ where $e_{Q, f}$ is the ramification index of $f$ at $Q$. 
Thus we get 
\[
\Delta_{Q, f^{\ast}L} = \Delta_{P, L} \cdot e_{Q, f}, 
\]
whence having 
\[
\sum_{Q \mapsto P} \Delta_{Q, f^{\ast}L} = \mathrm{deg}(f) \cdot \Delta_{P, L}. 
\]

Now let $S$ be finite subset of $\mathbb{P}^{1}$ and put 
\[
\Delta (L, S) = \sum_{P \in S} (\Delta_{P, L} - 1) . 
\]

\noindent
Here we have 

\begin{equation}
\Delta (f^{\ast}L, f^{-1}(S)) + \# f^{-1}(S) 
= \mathrm{deg}(f) \cdot \{ \Delta(L, S) + \# S \}, 
\label{eq:3-1-1}
\end{equation}

\noindent
and by Riemann-Hurwitz formula, 

\begin{equation}
-2 + 2 \cdot \mathrm{deg}(f) = \mathrm{deg}(f) \cdot (\# S) - \# f^{-1}(S). 
\label{eq:3-1-2}
\end{equation}

\noindent
When we take $\# S$ sufficiently large, (\ref{eq:3-1-1}) and (\ref{eq:3-1-2}) imply 
\[
\mathrm{deg}(f) \cdot (\Delta_{L} + 2) = \Delta_{f^{\ast}L} + 2, 
\]
and the lemma follows. 
}
}

\vspace{1ex}
Let $L_{n}$ denote a Lame operator with projective octahedral monodromy. 
By Theorem \ref{Litcanu}, we have $n \in \tfrac{1}{2}(\mathbb{N} + \tfrac{1}{2}) 
\cup \tfrac{1}{3}(\mathbb{N} + \tfrac{1}{2})$. 
As we saw in \S2.1, there is a $\ast$-morphism such that 
$L_{n} = f^{\ast}H_{1/2, 1/3, 1/4}$. 
We want to construct such a $\ast$-morphism $f : \mathbb{P}^{1} \rightarrow \mathbb{P}^{1}$.

In the octahedral case, by Lemma \ref{Bal-Dw}, 
\[
\mathrm{deg}(f) = 12n, 
\] 
and Riemann-Hurwitz formula implies 
\[
\# f^{-1} (\{0, 1, \infty \}) = 12 n + 2.
\] 
Then we can assume 
$f^{-1}(\{0, 1, \infty \}) = \{0, 1, \lambda, \infty, a_{1}, \cdots a_{12n-2} \}$ 
where, $a_{1}, \cdots a_{12n-2}$ denote distinct points different from $0, 1, \lambda, \infty$, 
and thus, possible ramification data of such an $f$ is given as follows: 

\begin{center}
  \begin{tabular}{|c|c|c|c|c|c|c|}                                                                       \hline
     {}      & 0    & 1    & $\lambda$ & $\infty$               & $a_{1}, \cdots a_{12n-2}$ & deg     \\ \hline 
     0       & 0, 1 & 0, 1 & 0, 1      & 0, $2n+1$              & 0, 2                      & $12n$   \\ \hline
     1       & 0    & 0    & 0         & 0, $3n+ \tfrac{3}{2}$  & 0, 3                      & $12n$   \\ \hline
     $\infty$& 0, 2 & 0, 2 & 0, 2      & 0, $4n+2$              & 0, 4                      & $12n$   \\ \hline
  \end{tabular}
\end{center}

\vspace{1ex}
Here, we explain how to read this table. 
For $P \in \{0, 1, \infty \}$ (an entry of the first column) and 
$Q \in \{0, 1, \infty, a_{1}, \cdots a_{12n-2} \}$ (an entry of the first row), 
the possible ramification index of $f$ at $P$ is written in the corresponding entry 
(i.e., $(Q, P)$-th entry); the number 0 occurs when $f(Q) \neq P$. 
These values $e_{Q, P}$ are calculated by the formula 
\[
\Delta_{Q, L_{n}} = e_{Q, f} \cdot \Delta_{P, H} 
\]
where $H = H_{1/2, 1/3, 1/4}$. 
Moreover, these values must satisfy the following compatibility conditions:

$(1)$ The summation of every row is equal to deg$(f)$. 

$(2)$ Every column contains only one non-zero number. 

\vspace{1ex}
Let us ask, conversely, if we can construct the $\ast$-morphism $f$, 
or what amounts to the same, the corresponding dessin, starting from a table 
as above which satisfies the above compatibility conditions for each
$n \in \tfrac{1}{2}( \mathbb{N} + \tfrac{1}{2}) \cup \tfrac{1}{3}( \mathbb{N} + \tfrac{1}{2})$.

\vspace{1ex}
\underline{\textbf{(1) The case for $n \in \tfrac{1}{2}( \mathbb{N} + \tfrac{1}{2}$)}}

\begin{center}
  \begin{tabular}{|c|c|c|c|c|c|c|}                                                                                       \hline
     {}      & 0 & 1 & $\lambda$ & $\infty$ & $a_{1}, \cdots a_{12n-2}$                 & deg     \\ \hline 
     0       & 1 & 1 & 1         & 0        & $( 6n-\tfrac{3}{2})$ pts with $mult. = 2$ & $12n$   \\ \hline
     1       & 0 & 0 & 0         & 0        & $4n$ pts with $mult. = 3$                 & $12n$   \\ \hline
     $\infty$& 0 & 0 & 0         & $4n+2$   & $( 2n-\tfrac{1}{2})$ pts with $mult. = 4$ & $12n$   \\ \hline
  \end{tabular}
\end{center}

If the dessin has $N$ loops with valency 4, 
\[
N = 2n - \tfrac{1}{2} \Leftrightarrow n = \tfrac{1}{2}(N + \tfrac{1}{2}) \in 
\tfrac{1}{2}(\mathbb{N} + \tfrac{1}{2}) 
\] 
and then the table becomes the following. 

\begin{center}
  \begin{tabular}{|c|c|c|c|c|c|c|}                                                                                       \hline
     {}      & 0 & 1 & $\lambda$ & $\infty$ & $a_{1}, \cdots a_{12n-2}$         &   deg        \\ \hline 
     0       & 1 & 1 & 1         & 0        & $3N$ pts with $mult. = 2$         & $ 6N + 3 $   \\ \hline
     1       & 0 & 0 & 0         & 0        & $( 2N + 1 )$ pts with $mult. = 3$ & $ 6N + 3 $   \\ \hline
     $\infty$& 0 & 0 & 0         & $2N+3$   & $  N  $ pts with $mult. = 4$      & $ 6N + 3 $   \\ \hline
  \end{tabular}
\end{center}
So it suffices to construct the dessins compatible with the table 
for all $N \in \mathbb{N}$. 
Now, we construct dessins.

\vspace{1ex} 
For $N = 0$,  
\begin{center}
\includegraphics[width=10em,clip]{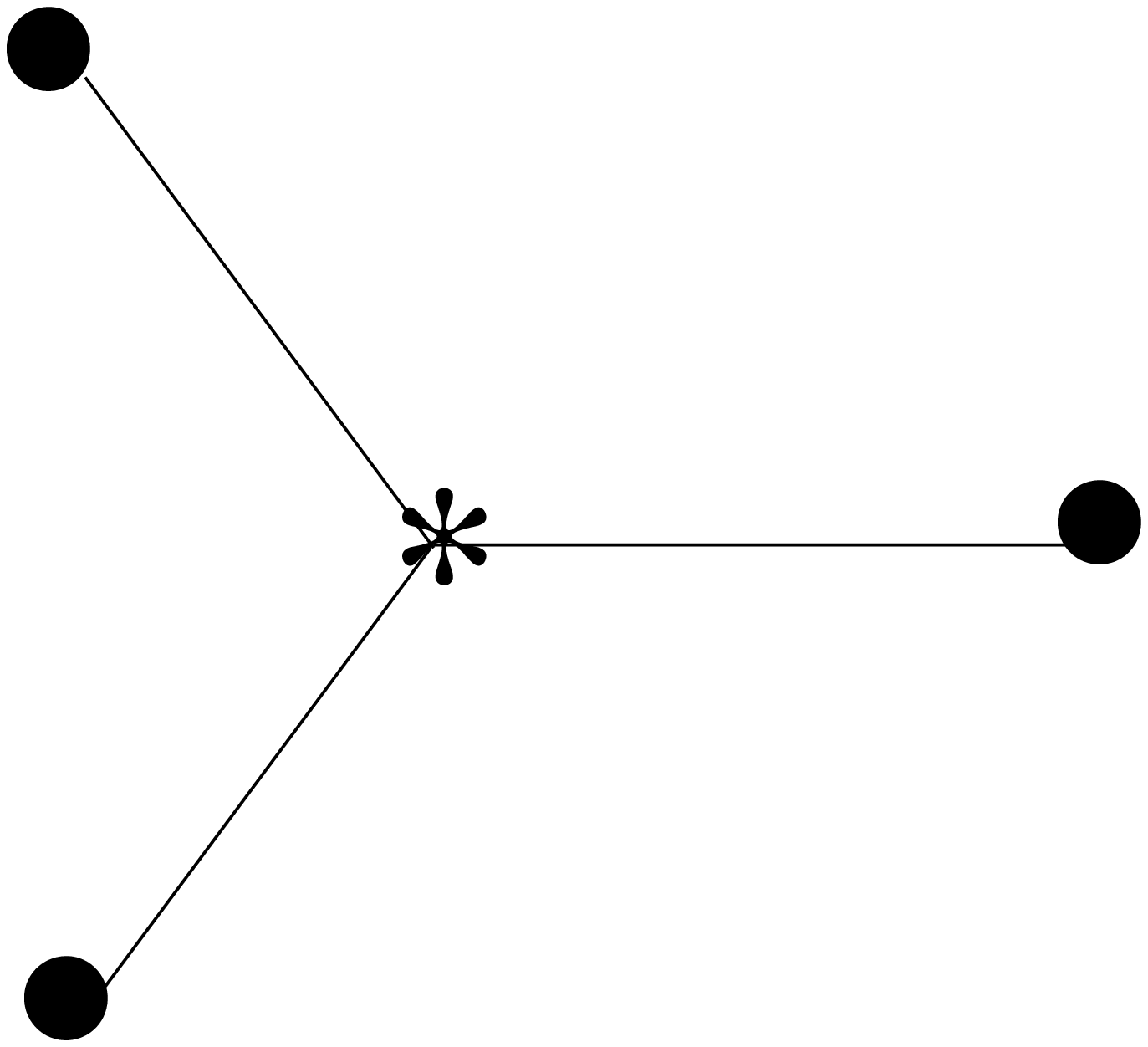}
\end{center}

For $N=1$, 
\begin{center}
\includegraphics[width=12em,clip]{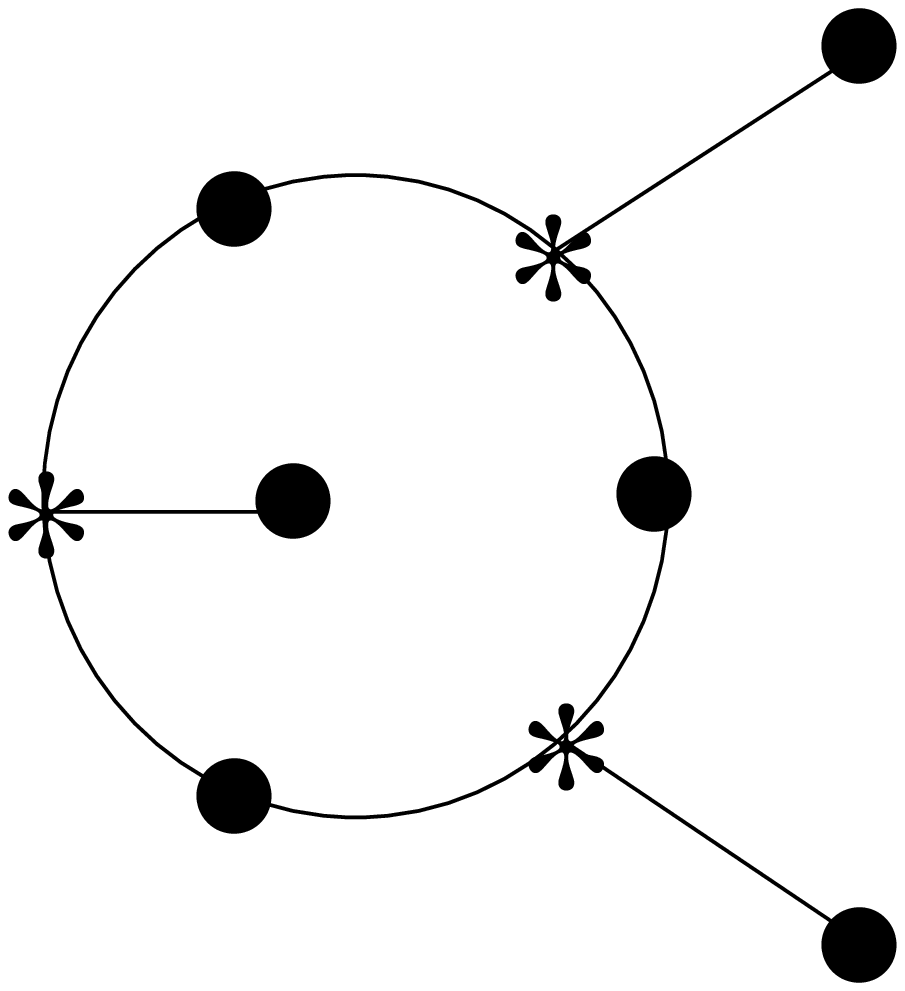}
\end{center}

For $N=2$,  
\begin{center}
\includegraphics[width=12em,clip]{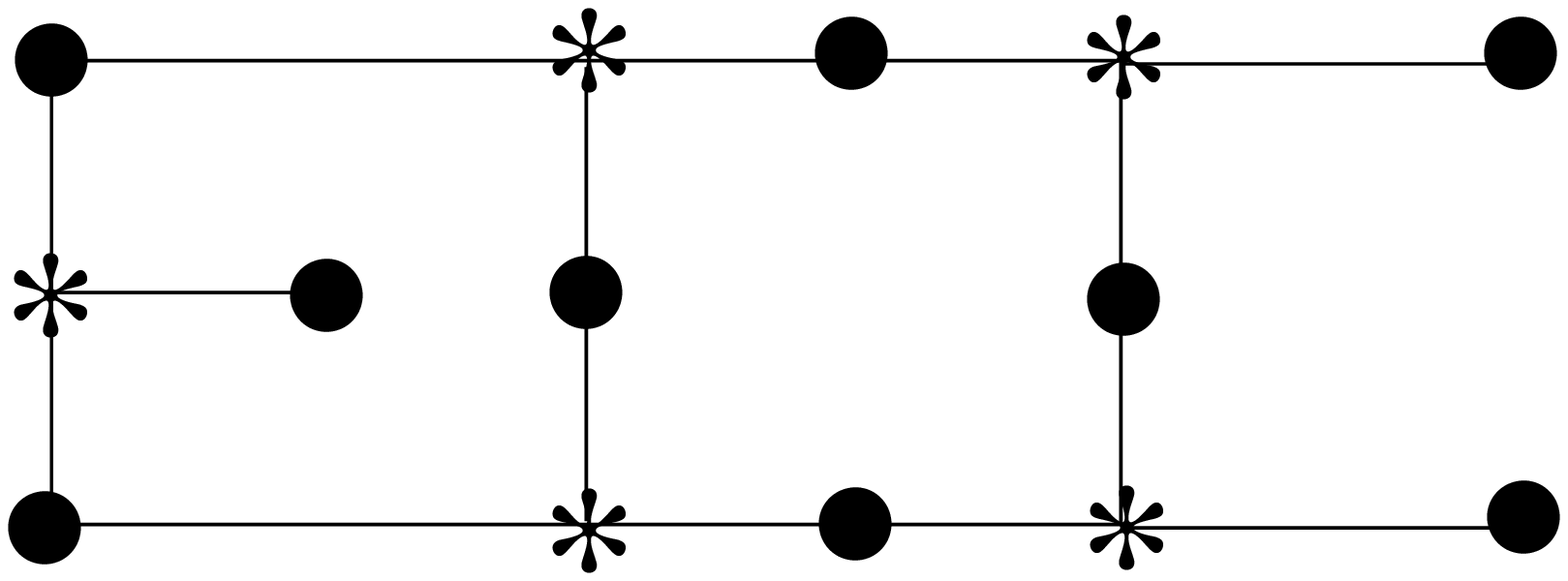}
\end{center}

For $N=k \geq 2$, 
\begin{center}
\includegraphics[width=15em,clip]{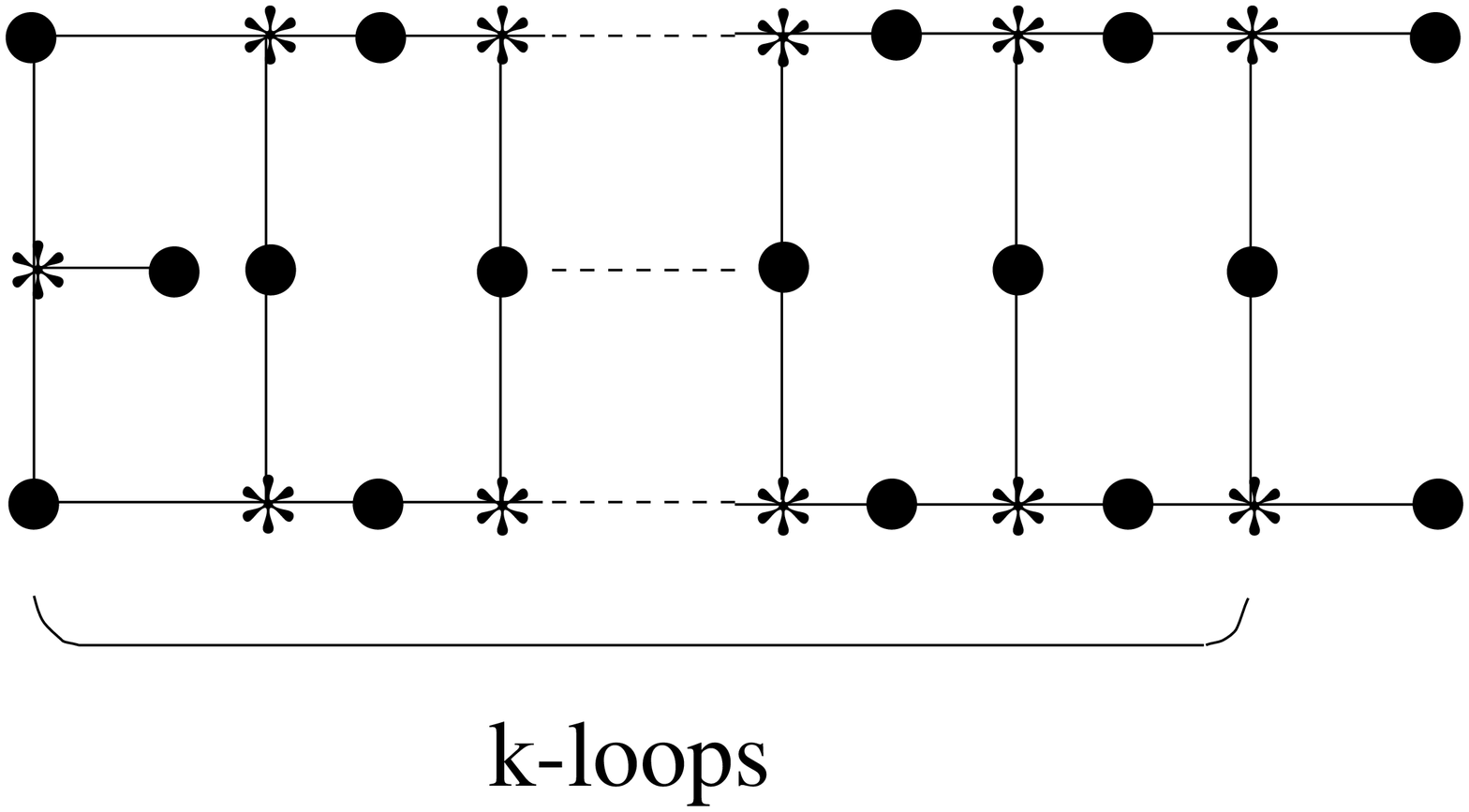}
\end{center}

It is easy to check that these are the dessins we want, and thus we could draw dessins inductively for all $N \in \mathbb{N}$. 

\begin{remark}
\rm{
For each $N \geq 1$, our dessin drawn above is one of those which are 
compatible with the table. 
There may exist other dessins compatible with the table. 
}
\end{remark}

\underline{\textbf{(2) The case for $n \in \tfrac{1}{3}( \mathbb{N} + \tfrac{1}{2}$)}} 

\begin{center}
  \begin{tabular}{|c|c|c|c|c|c|c|}                                                                            \hline
     {}      & 0 & 1 & $\lambda$ & $\infty$          & $a_{1}, \cdots a_{12n-2}$                 & deg     \\ \hline 
     0       & 1 & 1 & 0         & 0                 & $( 6n-1 )$ pts with $mult. = 2$           & $12n$   \\ \hline
     1       & 0 & 0 & 0         & $3n+\tfrac{3}{2}$ & $( 3n-\tfrac{1}{2})$ pts with $mult. = 3$ & $12n$   \\ \hline
     $\infty$& 0 & 0 & 2         & 0                 & $( 3n-\tfrac{1}{2})$ pts with $mult. = 4$ & $12n$   \\ \hline
  \end{tabular}
\end{center}

If the dessin has $N$ loops with valency 4, 
\[
N = 3n - \tfrac{1}{2} \Leftrightarrow n = \tfrac{1}{3}(N + \tfrac{1}{2}) \in 
\tfrac{1}{3}(\mathbb{N} + \tfrac{1}{2}) 
\] 
and then the table becomes the following. 

\begin{center}
  \begin{tabular}{|c|c|c|c|c|c|c|}                                                                            \hline
     {}      & 0 & 1 & $\lambda$ & $\infty$ & $a_{1}, \cdots a_{12n-2}$ & deg      \\ \hline 
     0       & 1 & 1 & 0         & 0        & $2N$ pts with $mult. = 2$ & $4N + 2$ \\ \hline
     1       & 0 & 0 & 0         & $N + 2$  & $N$ pts with $mult. = 3$  & $4N + 2$ \\ \hline
     $\infty$& 0 & 0 & 2         & 0        & $N$ pts with $mult. = 4$  & $4N + 2$ \\ \hline
  \end{tabular}
\end{center}
As in the previous case, it suffices to construct the dessins compatible with the table 
for all $N \in \mathbb{N}$. 

Now, we construct dessins. 
(For each $N \geq 1$, our dessin is one of those which are compatible with the table.) 

\vspace{1ex}
For $N = 0$, 
\begin{center}
\includegraphics[width=10em,clip]{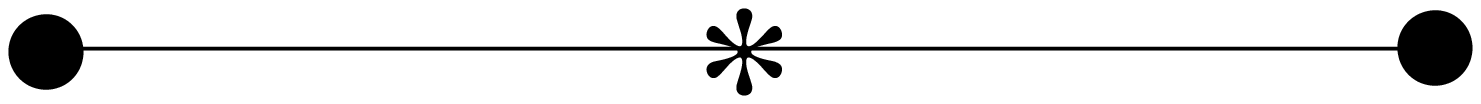}
\end{center}

For $N=1$, 
\begin{center}
\includegraphics[width=12em,clip]{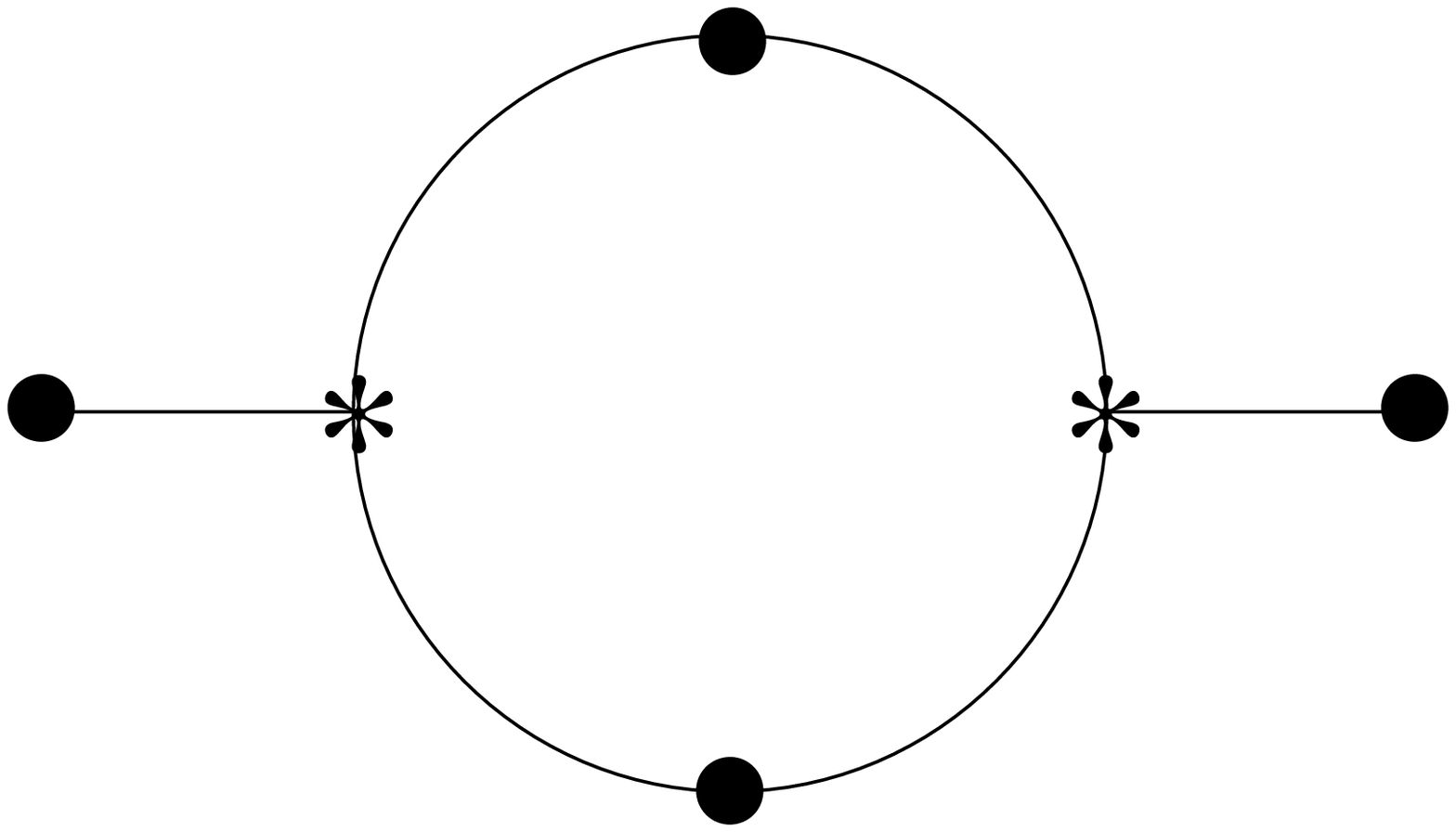}
\end{center}

For $N=2$, 
\begin{center}
\includegraphics[width=12em,clip]{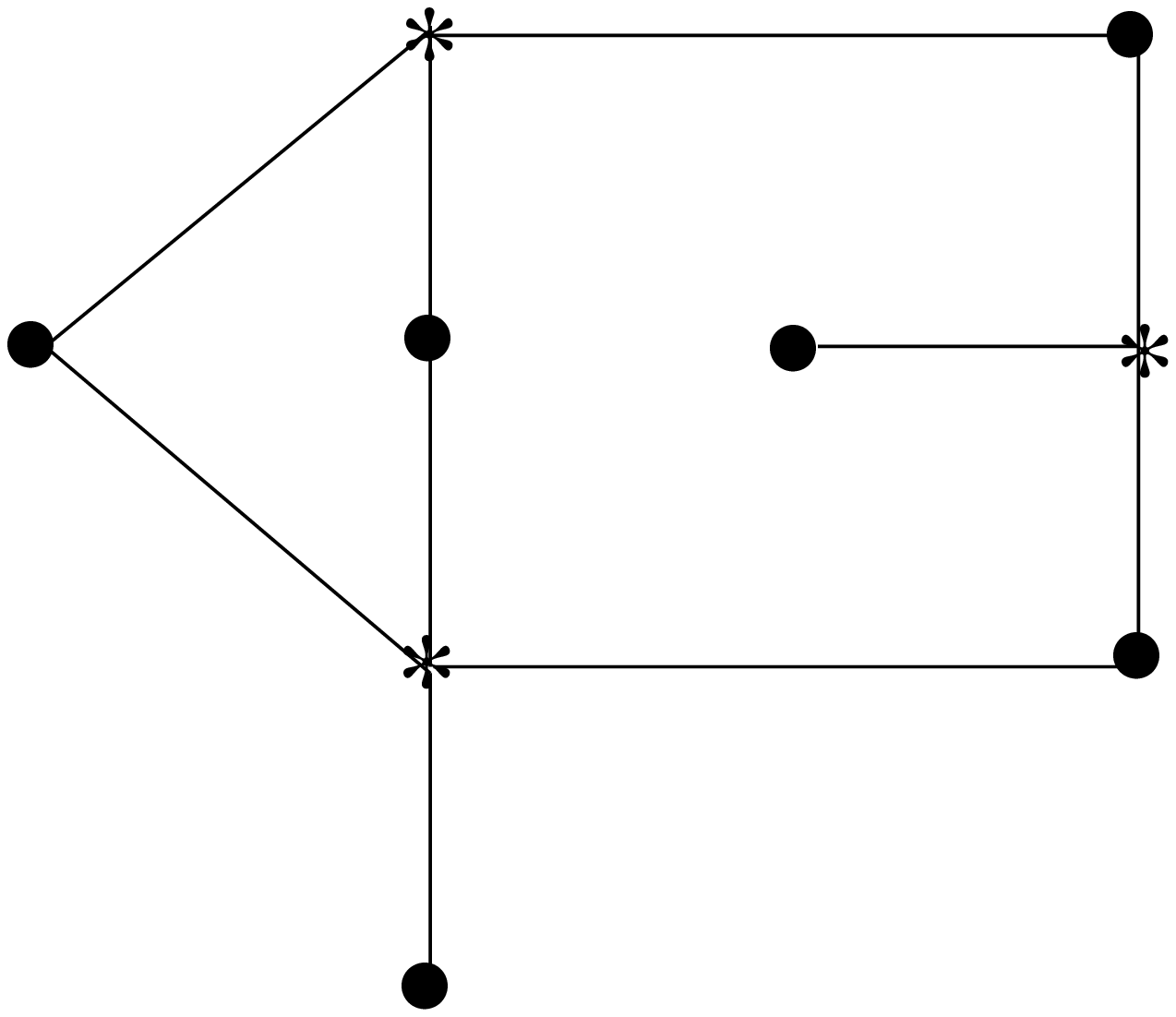}
\end{center}

For $N=3$, 
\begin{center}
\includegraphics[width=10em,clip]{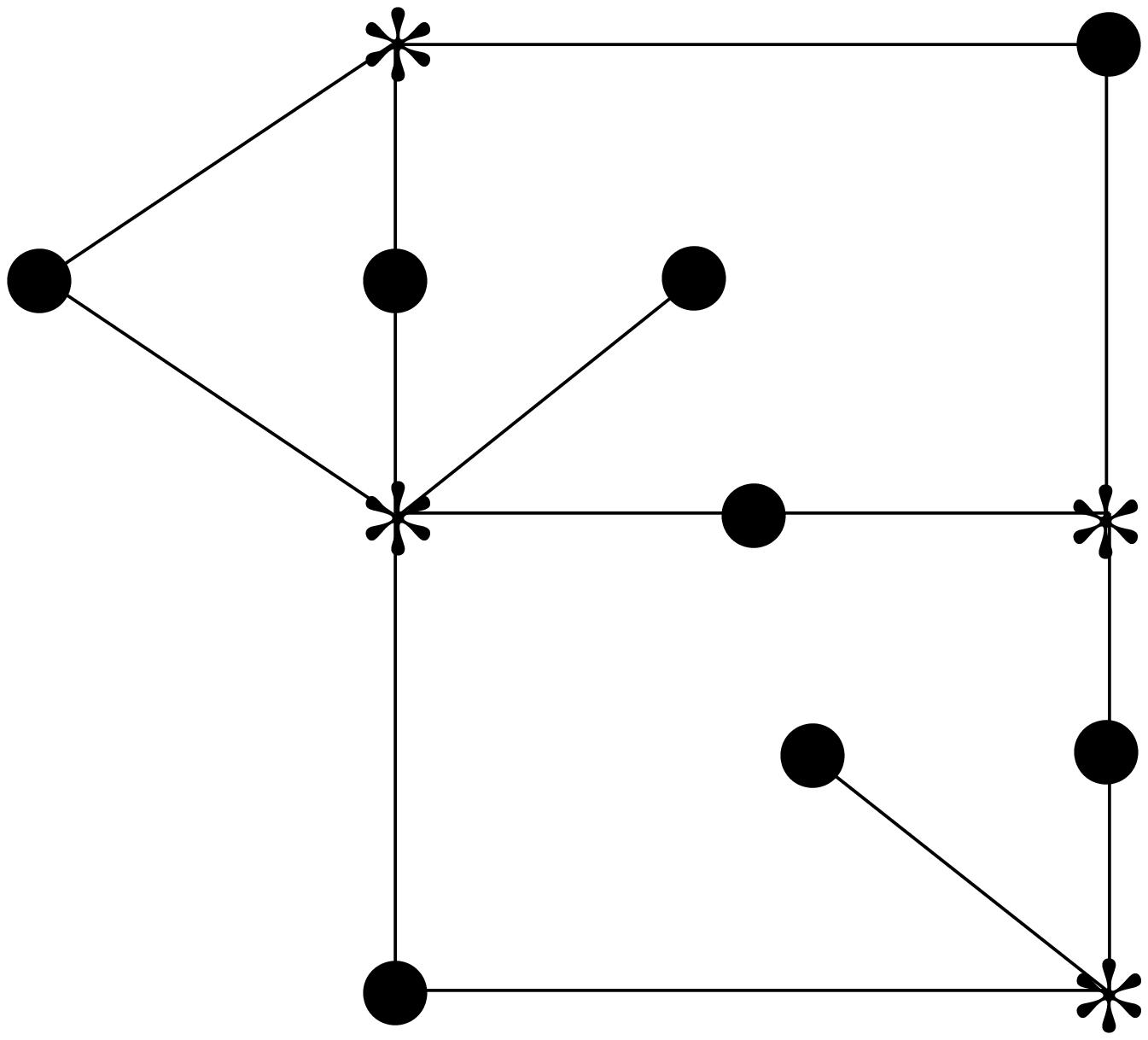}
\end{center}

For $N=4$, 
\begin{center}
\includegraphics[width=12em,clip]{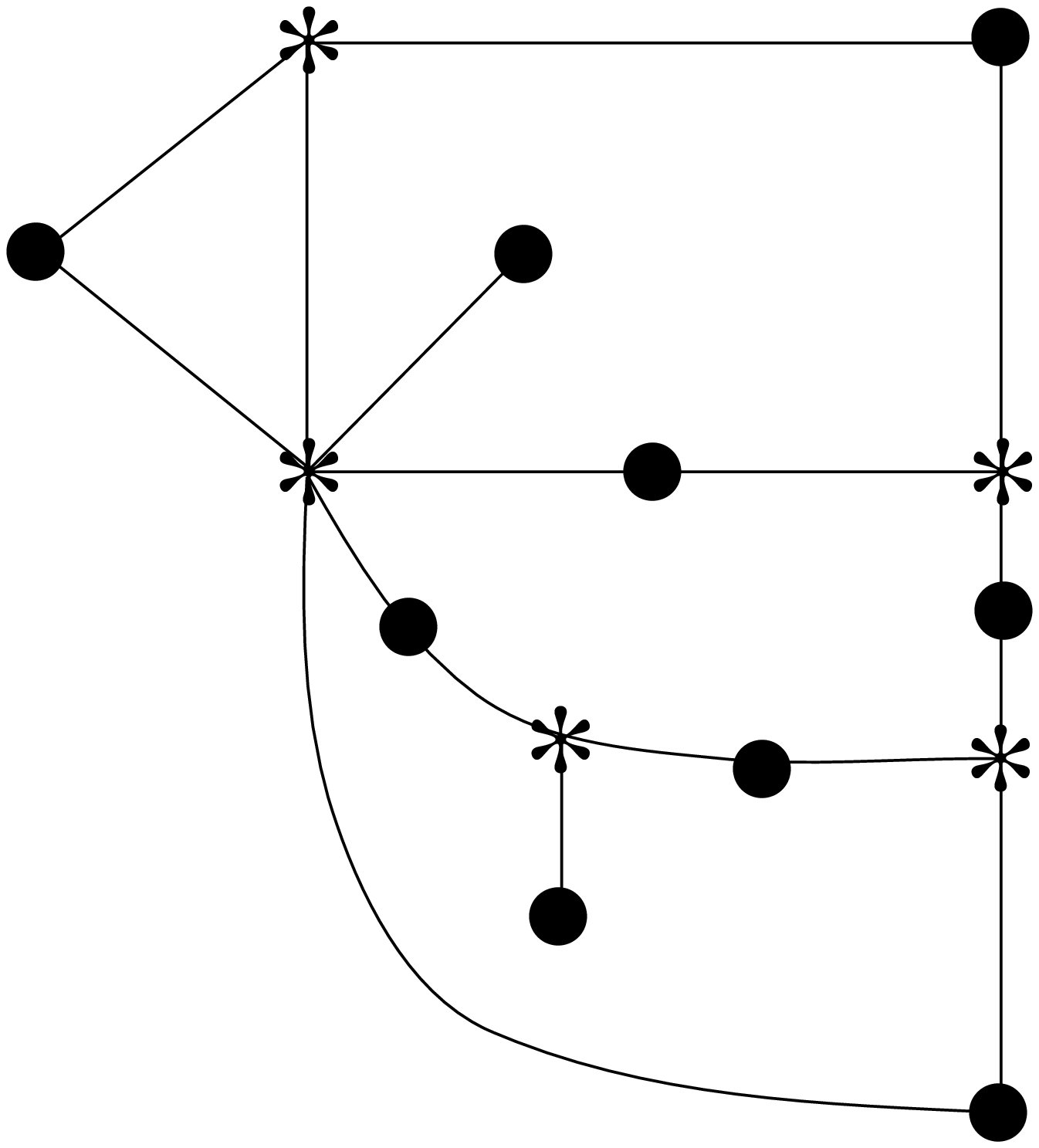}
\end{center}

For $N=5$, 
\begin{center}
\includegraphics[width=12em,clip]{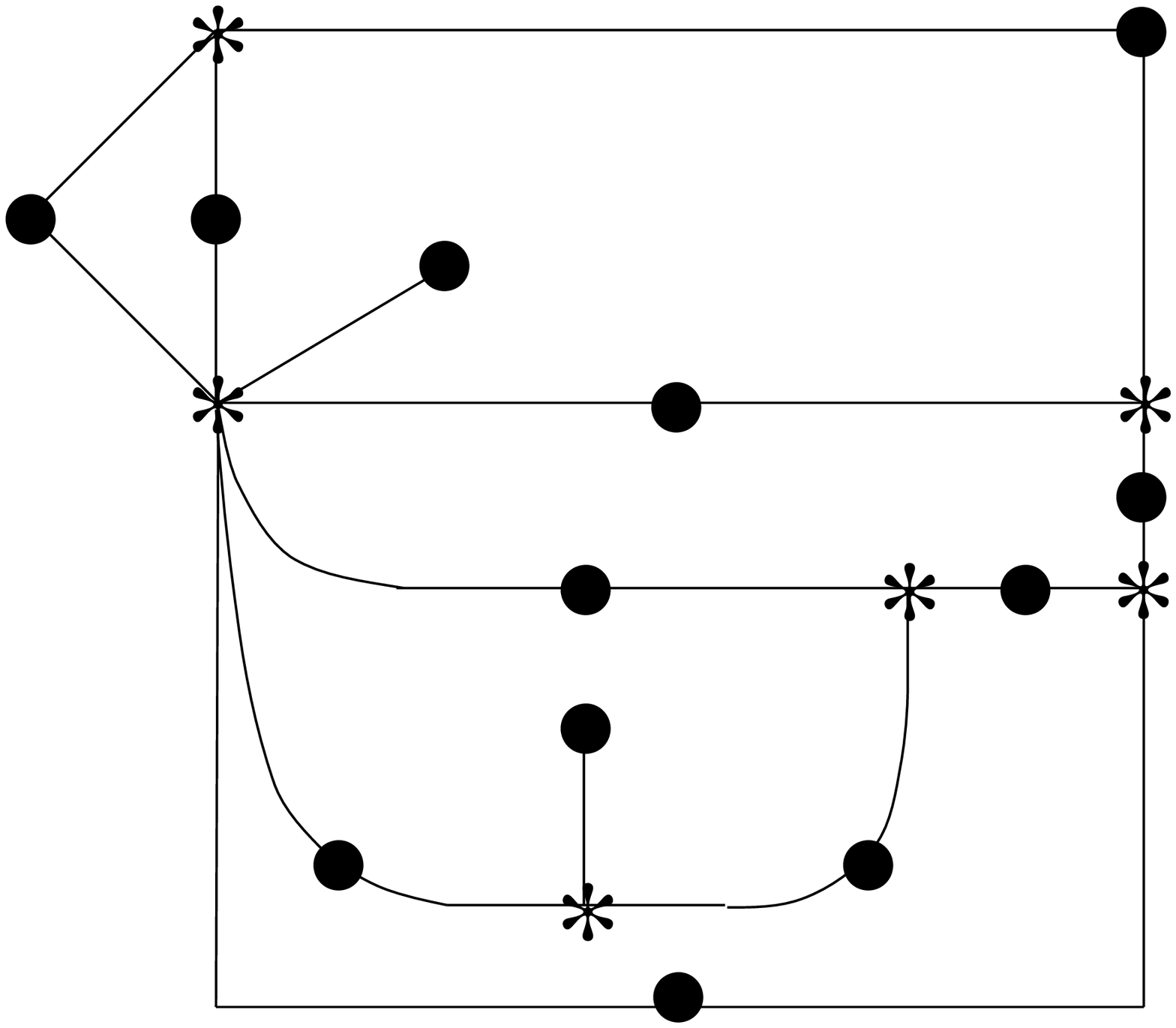}
\end{center}

For $N=6$, 
\begin{center}
\includegraphics[width=12em,clip]{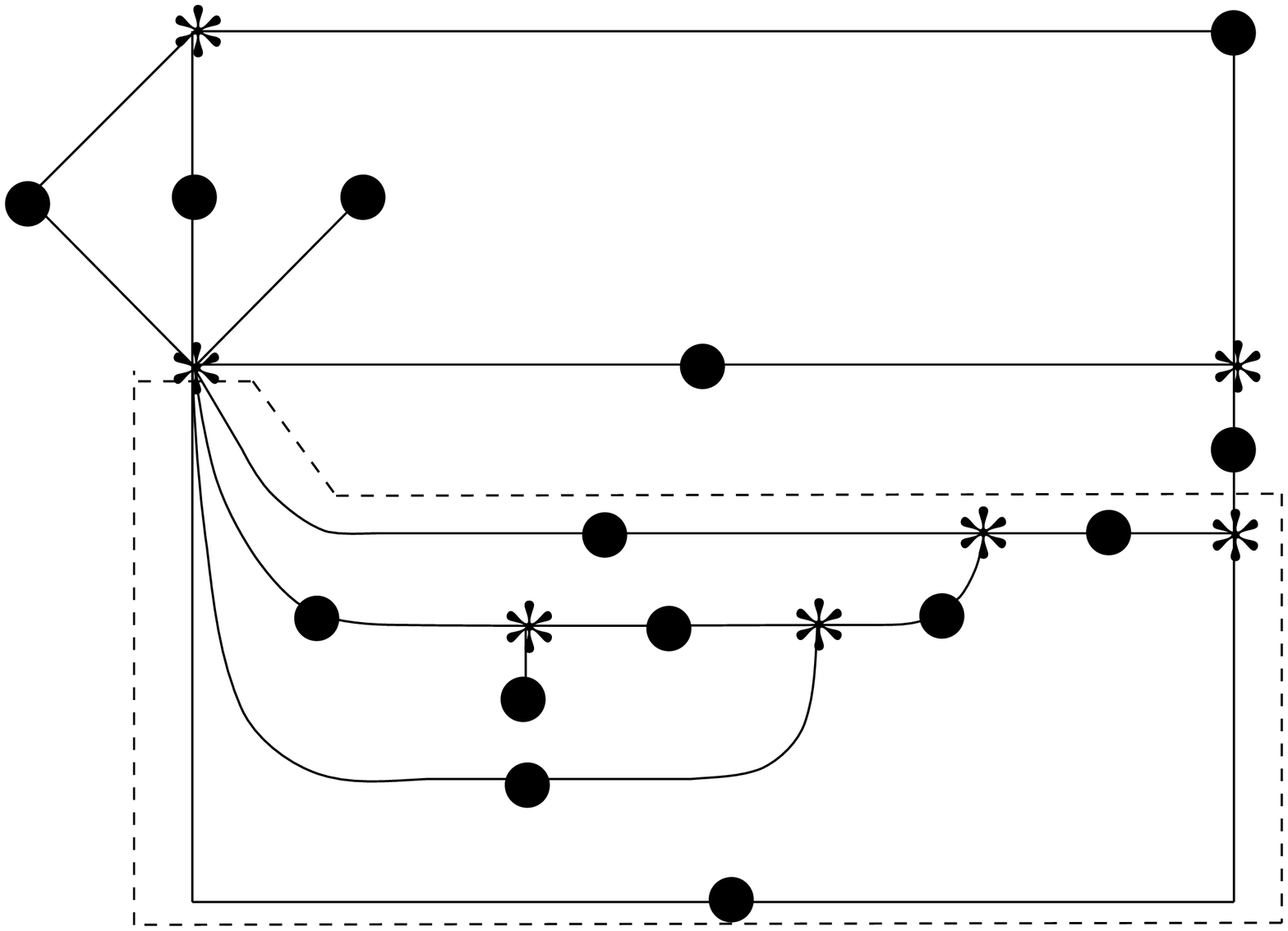}
\end{center}

For $N \geq 6$, we can construct dessins inductively according to 
the following operations. We operate the lower-half part of the dessin. 
(We draw only the part which is enclosed by the dotted line in the dessin of $N=6$.) 

If the dessin with $N$ loops is 
 
\begin{center}
\includegraphics[width=15em,clip]{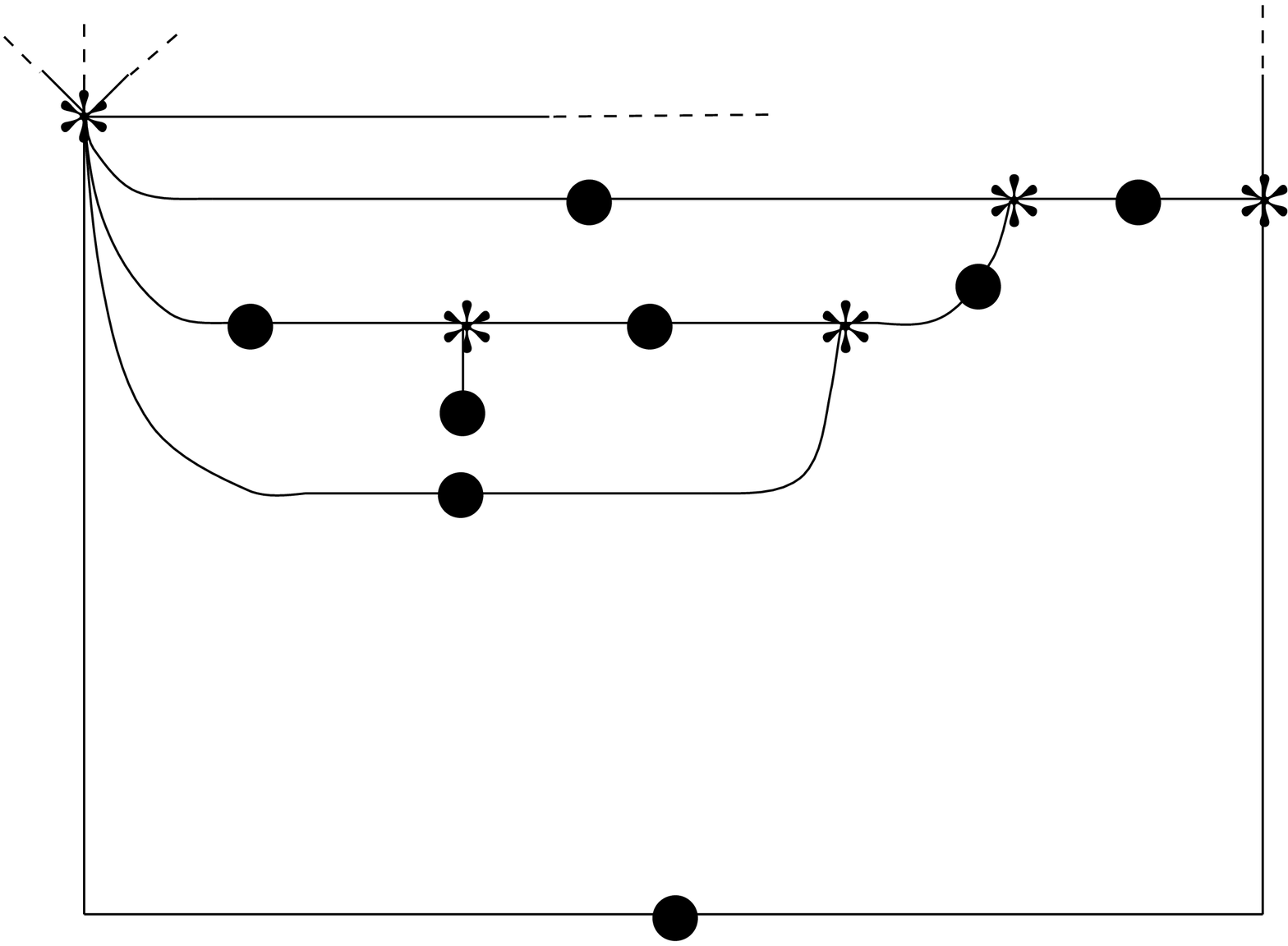}
\end{center}

\noindent
then the dessin with $N+1$ loops is obtained by  

\begin{center}
\includegraphics[width=15em,clip]{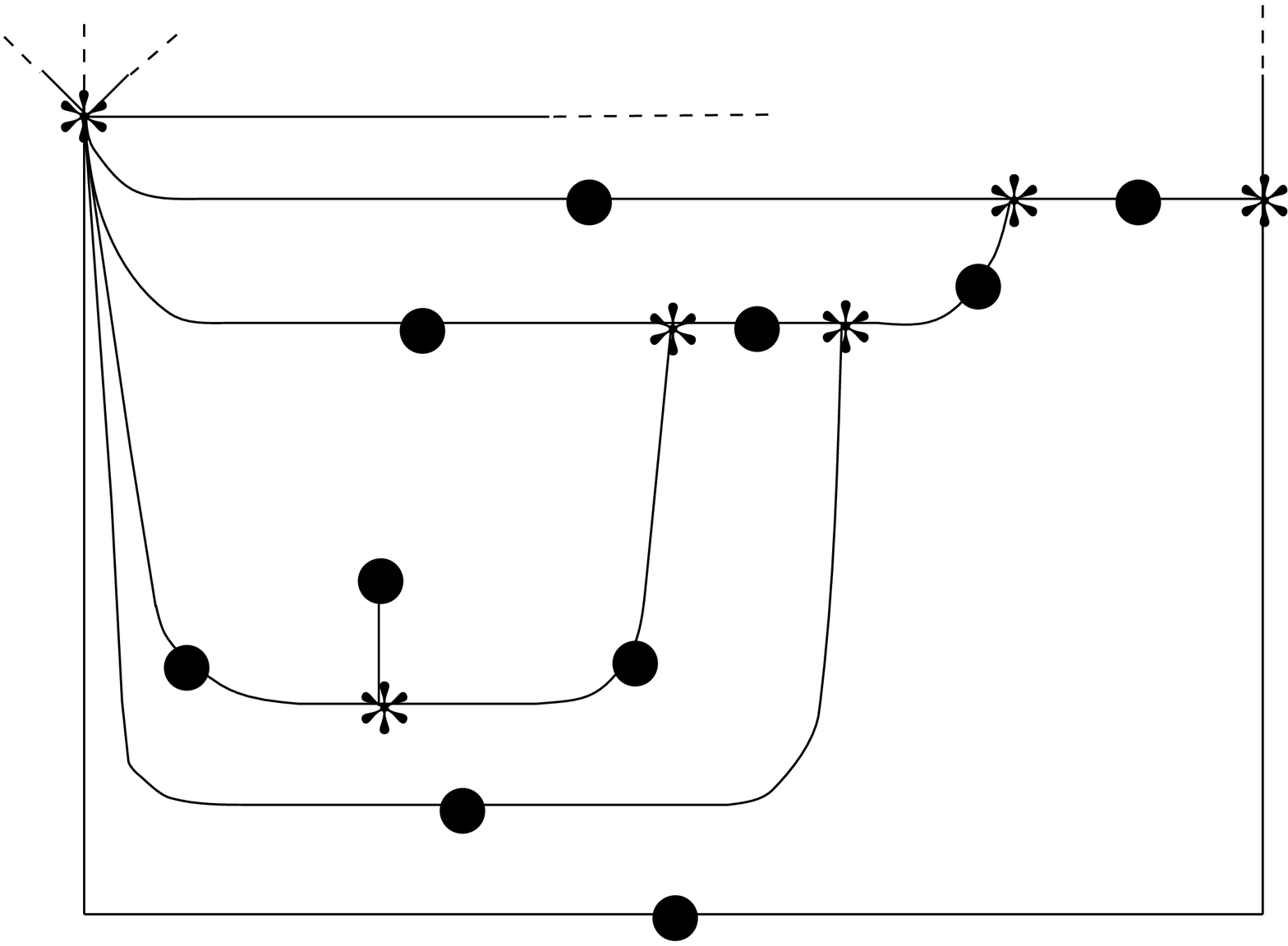}
\end{center}

\noindent
and the dessin with $N+2$ loops is obtained by 
 
\begin{center}
\includegraphics[width=15em,clip]{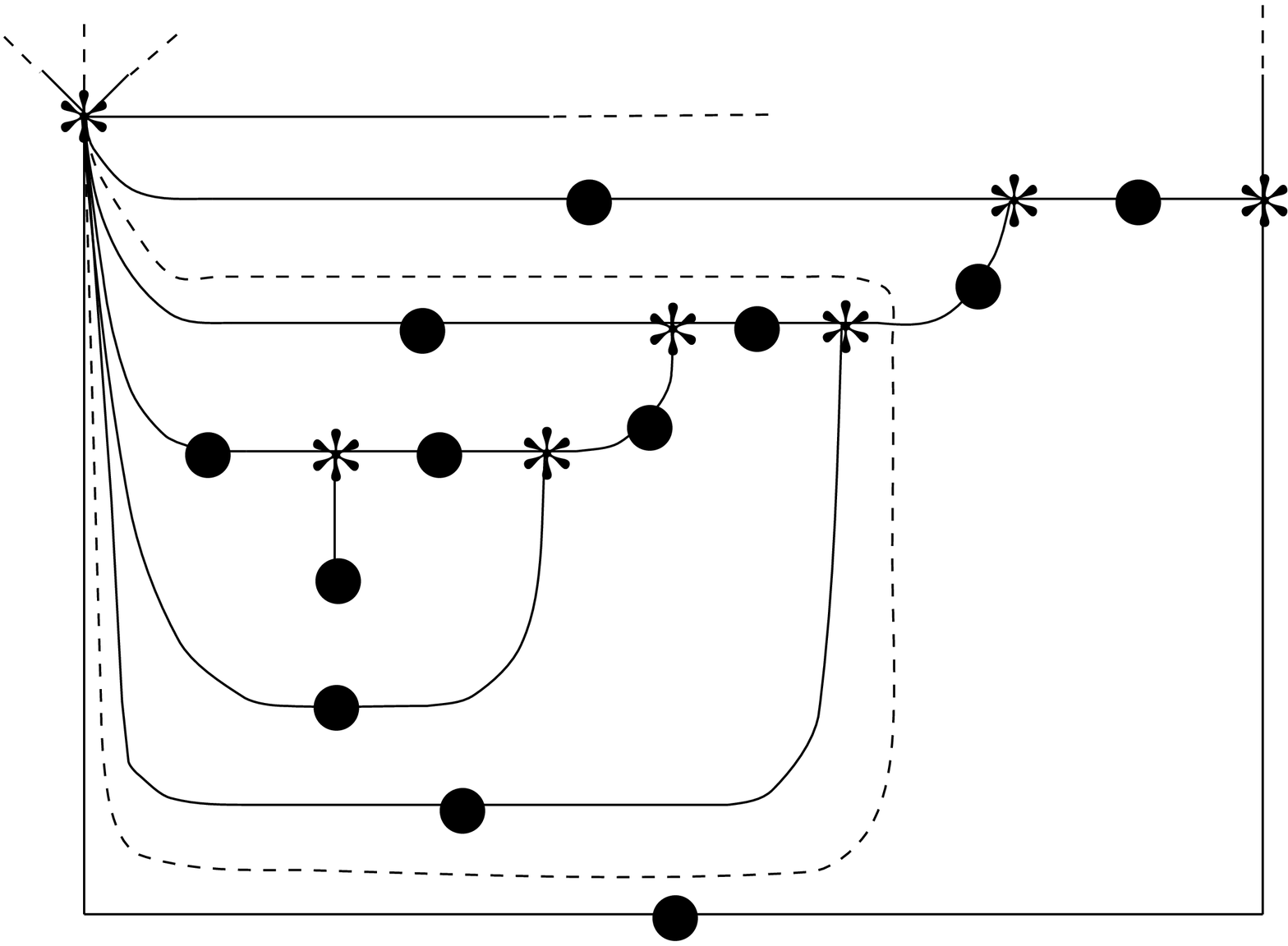}
\end{center}

We can see that the part which is enclosed by the dotted line 
in the above dessin repeatedly appears when $N$ is even. 
Also we can easily see that dessins inductively constructed by this operation 
are compatible with the table above. 
In fact, in each step of this operation, 
the number of $``\bullet"$ with valency 2 increase by two, 
that of $``\ast"$ with valency 3 increase by one, 
that of edges increase by four, and 
that of loops having $a_{i}$ in the fiber over $\infty$ with valency 4 inside 
increase by one. 

\begin{remark}\label{rem-exc1}
\rm{
As we saw in the introduction, we can assume degree parameter $n > -1/2$. 
If there exists Lam\'{e} operator with projective octahedral monodromy 
for $n=-1/4$ (resp. $n=-1/6$), 
there would exist Belyi morphism which pull-backs $H_{1/2, 1/3, 1/4}$ 
into $L_{-1/4}$ (resp. $L_{-1/6}$). 
But this morphism has to have negative degree, and hence the corresponding dessin does not exist. 
Therefore, the Lam\'{e} operator with projective octahedral monodromy 
for $n=-1/4$ and $n=-1/6$ do not exist. 
}
\end{remark}

\subsection{The case of projective icosahedral monodromy.}

As in the previous subsection, 
We want $\ast$ - morphism $f : \mathbb{P}^{1} \rightarrow \mathbb{P}^{1}$ such that 
$L_{n} = f^{\ast}H_{1/2, 1/3, 1/5}$ for 
$n \in \tfrac{1}{3}( \mathbb{N} + \tfrac{1}{2}) \cup \tfrac{1}{5}( \mathbb{N} + \tfrac{1}{2}) $. 

Lemma \ref{Bal-Dw} and Riemann-Hurwitz formula implies 
\[
\# f^{-1} (\{0, 1, \infty \}) = 30 n + 2.
\] 
So we can assume 
$f^{-1}(\{0, 1, \infty \}) = \{0, 1, \lambda, \infty, a_{1}, \cdots a_{30n-2} \}$, 
and possible ramification data of $f$ is according to the following table. 

\begin{center}
  \begin{tabular}{|c|c|c|c|c|c|c|}                                                                       \hline
     {}      & 0    & 1    & $\lambda$ & $\infty$               & $a_{1}, \cdots a_{12n-2}$ & deg     \\ \hline 
     0       & 1    & 1    & 1         & 0, $2n+1$              & 0, 2                      & $30n$   \\ \hline
     1       & 0    & 0    & 0         & 0, $3n+ \tfrac{3}{2}$  & 0, 3                      & $30n$   \\ \hline
     $\infty$& 0    & 0    & 0         & 0, $5n+ \tfrac{5}{2}$  & 0, 5                      & $30n$   \\ \hline
  \end{tabular}
\end{center}

In this subsection, we construct dessins compatible with the ramification data above for each 
$n \in \tfrac{1}{3}( \mathbb{N} + \tfrac{1}{2}) \cup \tfrac{1}{5}( \mathbb{N} + \tfrac{1}{2}) $. 

\vspace{1ex}
\underline{\textbf{(3) The case for $n \in \tfrac{1}{3}( \mathbb{N} + \tfrac{1}{2}$)}}

\begin{center}
  \begin{tabular}{|c|c|c|c|c|c|c|}                                                                                       \hline
     {}      & 0 & 1 & $\lambda$ & $\infty$          & $a_{1}, \cdots a_{12n-2}$                  & deg     \\ \hline 
     0       & 1 & 1 & 1         & 0                 & $( 15n-\tfrac{3}{2})$ pts with $mult. = 2$ & $30n$   \\ \hline
     1       & 0 & 0 & 0         & $3n+\tfrac{3}{2}$ & $( 9n-\tfrac{1}{2})$ pts with $mult. = 3$  & $30n$   \\ \hline
     $\infty$& 0 & 0 & 0         & 0                 & $6n$ pts with $mult. = 5$                  & $30n$   \\ \hline
  \end{tabular}
\end{center}

When we have $n \in \tfrac{1}{3}(\mathbb{N} + \tfrac{1}{2})$, 
$6n$ is odd and we can put $6n = 2M + 1$ where $M \in \mathbb{N}$. 
Then the table becomes the following. 

\begin{center}
  \begin{tabular}{|c|c|c|c|c|c|c|}                                                                                       \hline
     {}      & 0 & 1 & $\lambda$ & $\infty$ & $a_{1}, \cdots a_{12n-2}$         & deg        \\ \hline 
     0       & 1 & 1 & 1         & 0        & $( 5M + 1 )$ pts with $mult. = 2$ & $10 M + 5$ \\ \hline
     1       & 0 & 0 & 0         & $M + 2$  & $( 3M + 1 )$ pts with $mult. = 3$ & $10 M + 5$ \\ \hline
     $\infty$& 0 & 0 & 0         & 0        & $( 2M + 1 )$ pts with $mult. = 5$ & $10 M + 5$ \\ \hline
  \end{tabular}
\end{center}
So it suffices to construct the dessins compatible with the table 
for all $M \in \mathbb{N}$. 
Now, we construct dessins. 
(For each $M \geq 1$, our dessin is one of those which are compatible with the table.) 

\vspace{1ex}
For $M=0$, 
\begin{center}
\includegraphics[width=12em,clip]{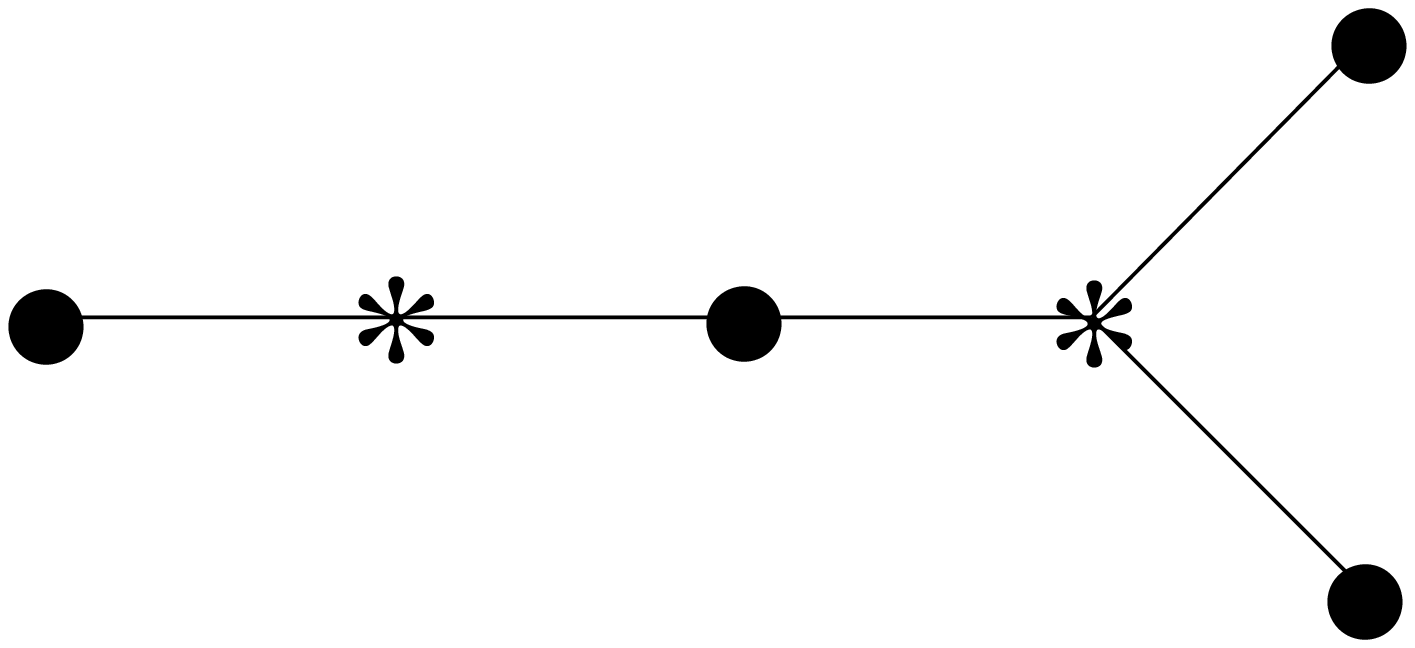}
\end{center}

For $M=1$, 
\begin{center}
\includegraphics[width=12em,clip]{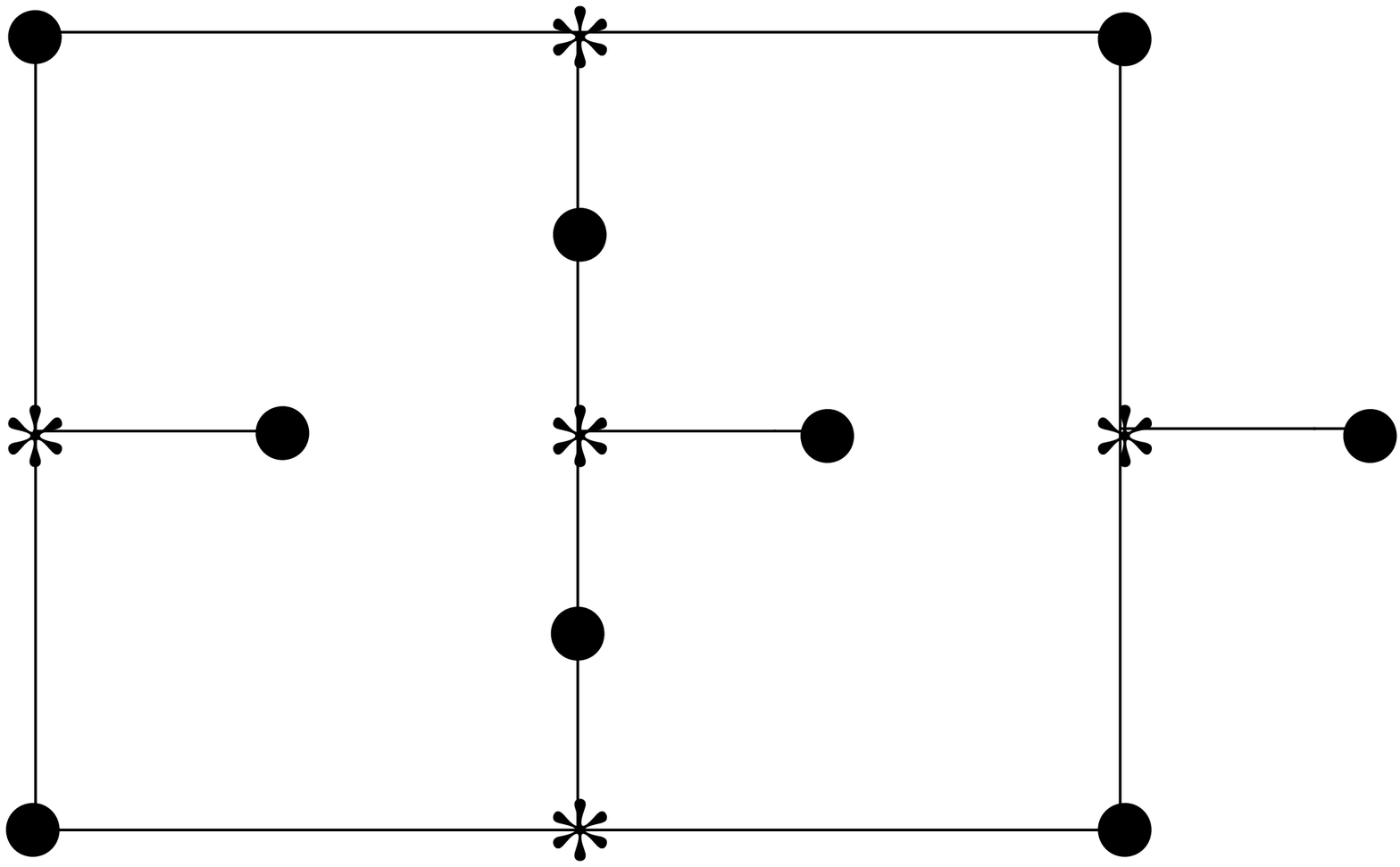}
\end{center}

For $M=2$, 
\begin{center}
\includegraphics[width=12em,clip]{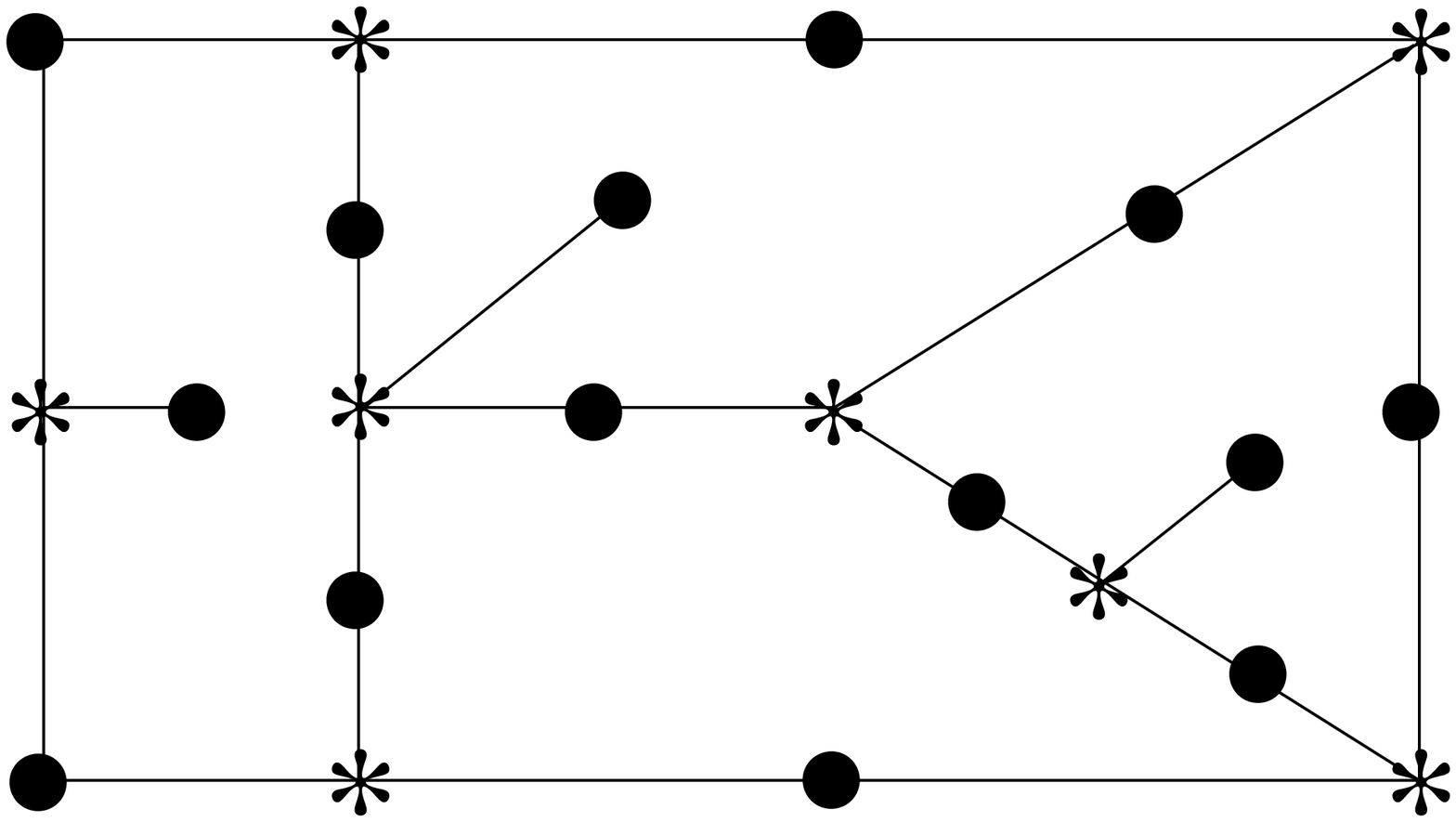}
\end{center}

For $M=3$, 
\begin{center}
\includegraphics[width=12em,clip]{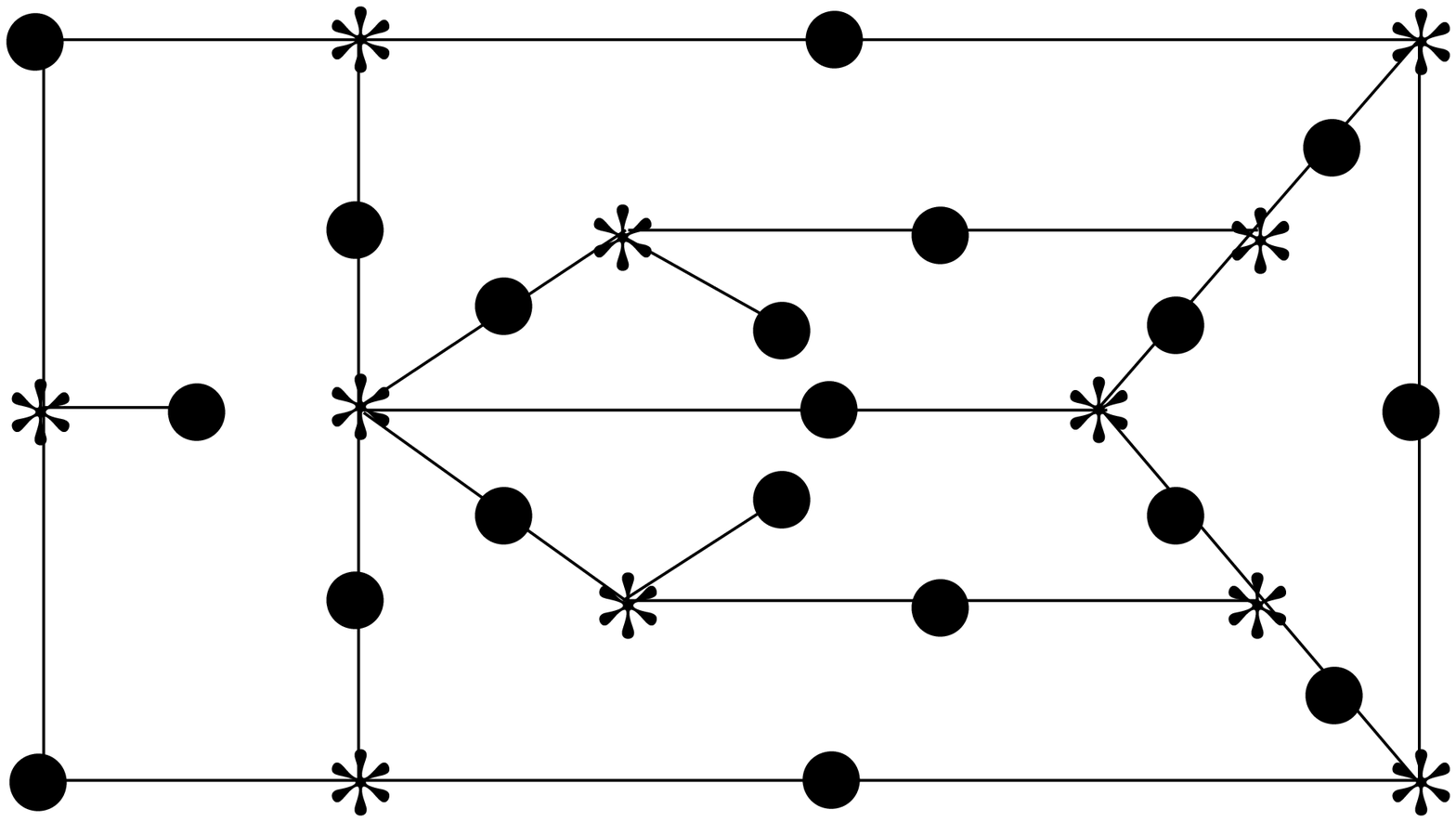}
\end{center}

For $M=4$, 
\begin{center}
\includegraphics[width=12em,clip]{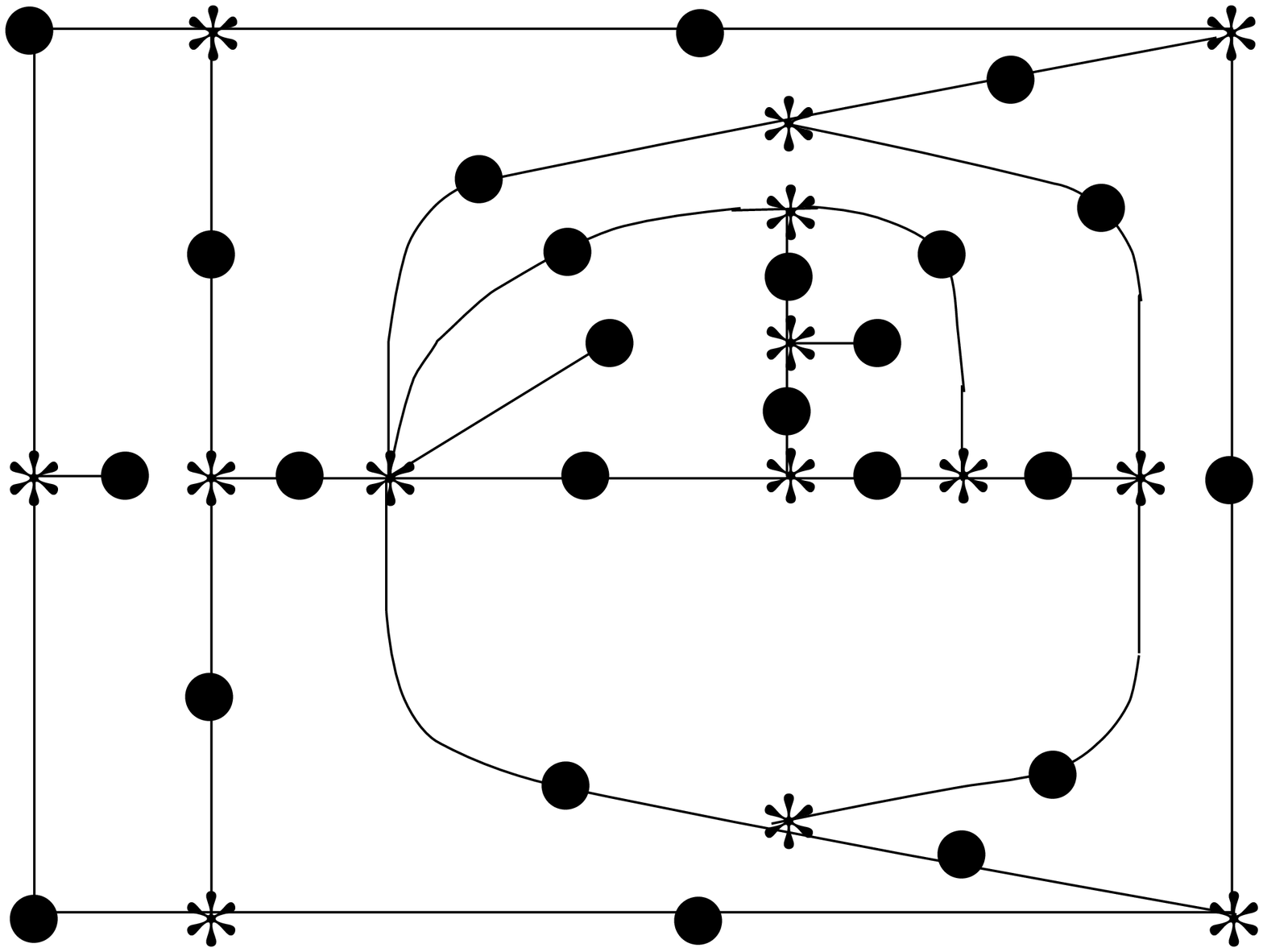}
\end{center}

For $M=5$, 
\begin{center}
\includegraphics[width=15em,clip]{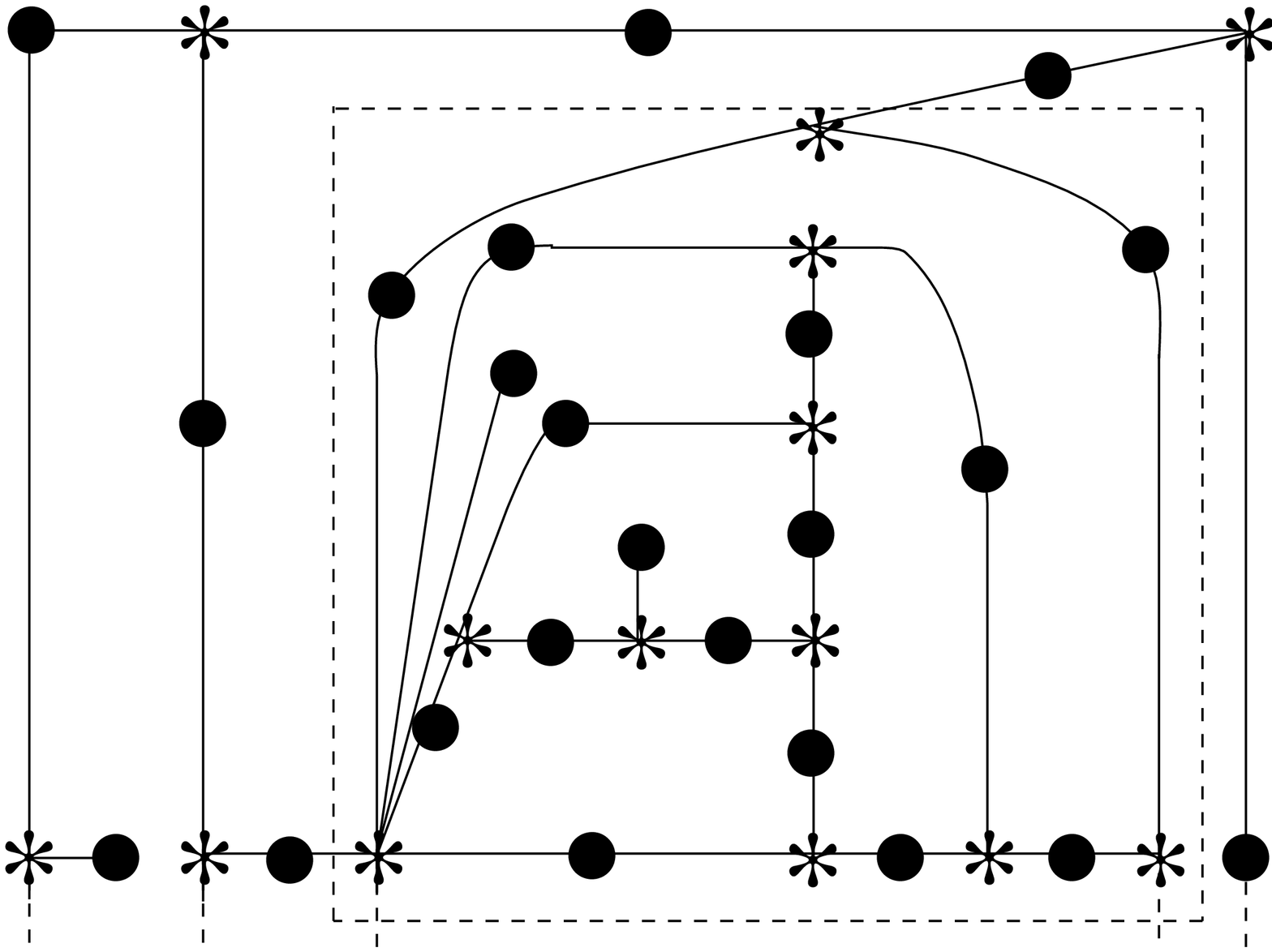}
\end{center}
(We draw only the upper-half part, the lower part looks like the same as $M=4$.) 

For $M=6$, 
\begin{center}
\includegraphics[width=15em,clip]{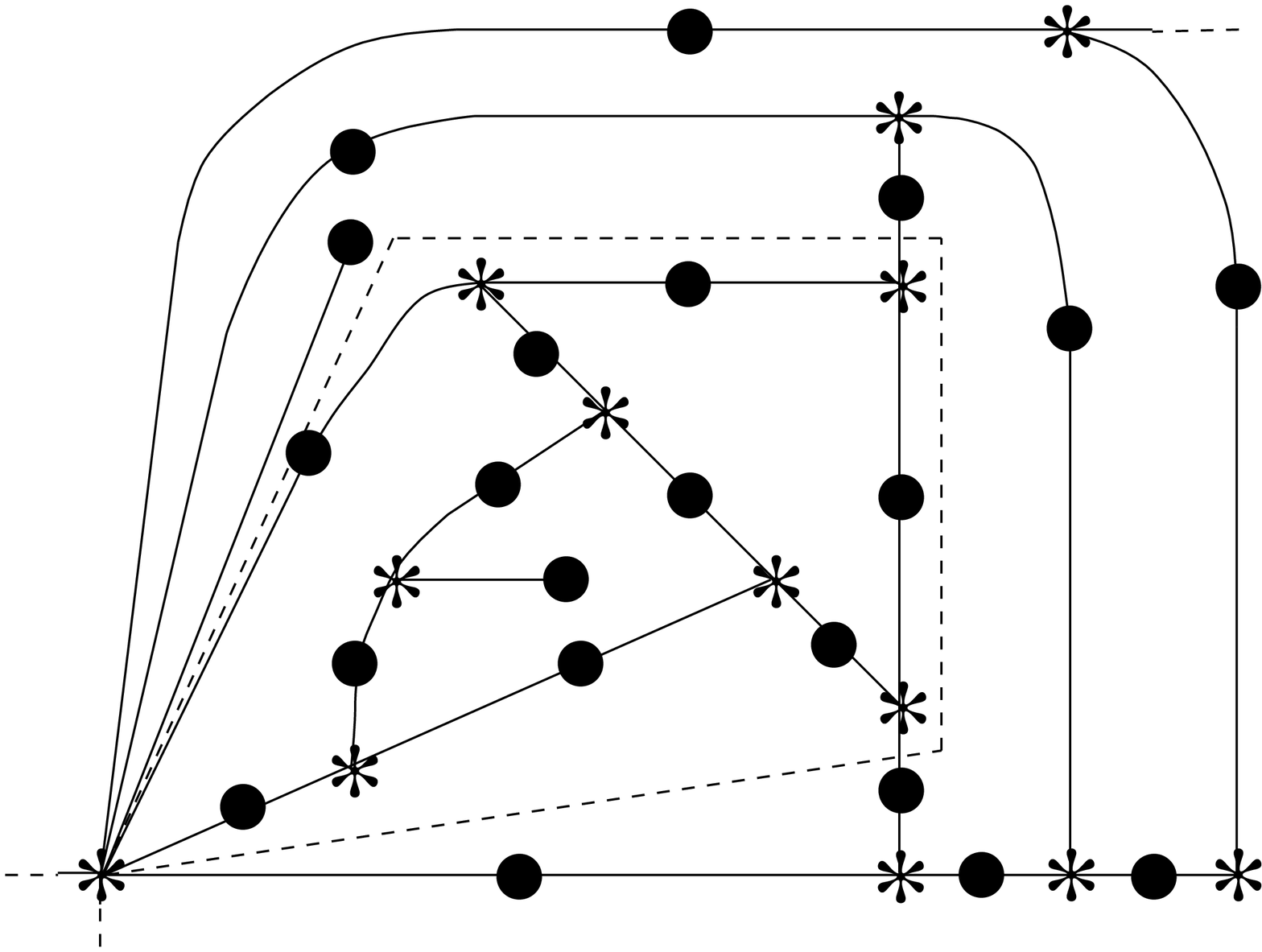}
\end{center}
(We draw only the main part which is enclosed by the dotted line in the dessin of $M=5$, 
the other parts are same as $M=5$.) 

\vspace{1ex}
For $M \geq 6$, we can construct dessins inductively according to 
the following operations. We operate the upper half part of the dessin. 
(We draw only the main part that is enclosed by the dotted line in the dessin of $M=6$.) 

If the dessin with $M=2m$ loops is 
 
\begin{center}
\includegraphics[width=15em,clip]{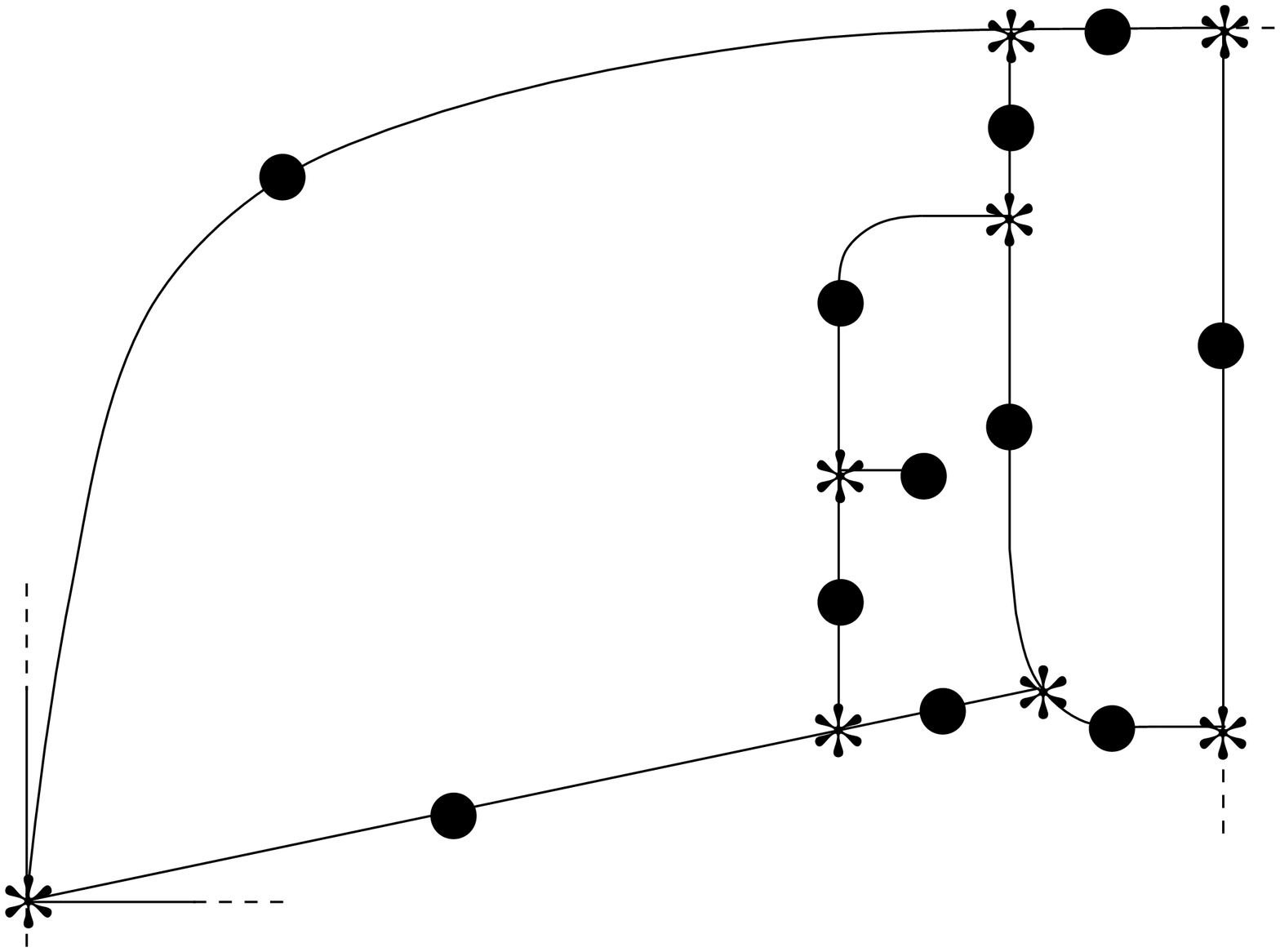}
\end{center}

\noindent
then the dessin with $2m+1$ loops is obtained by  

\begin{center}
\includegraphics[width=15em,clip]{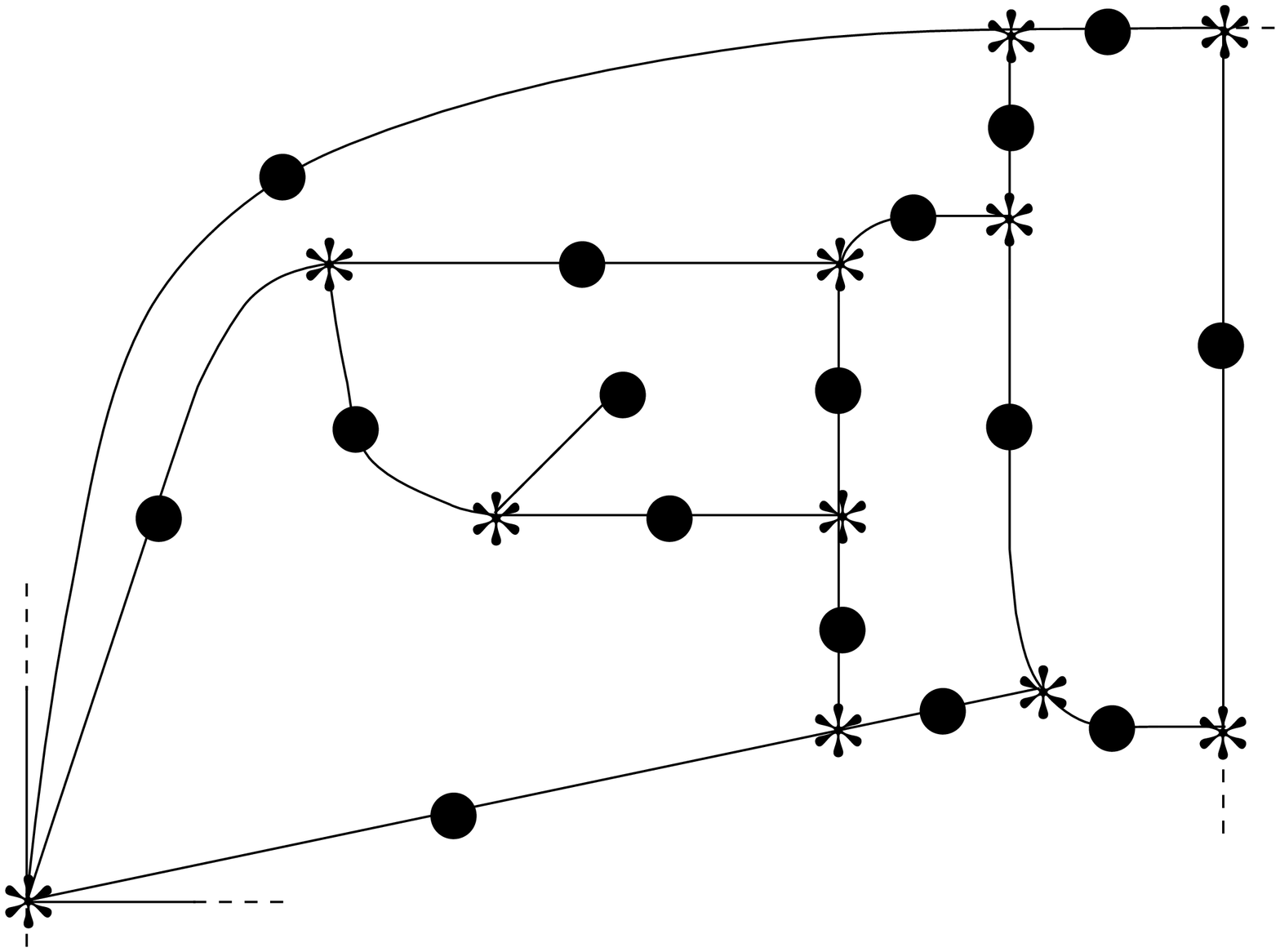}
\end{center}

\noindent
and the dessin with $2m+2$ loops is obtained by 
 
\begin{center}
\includegraphics[width=15em,clip]{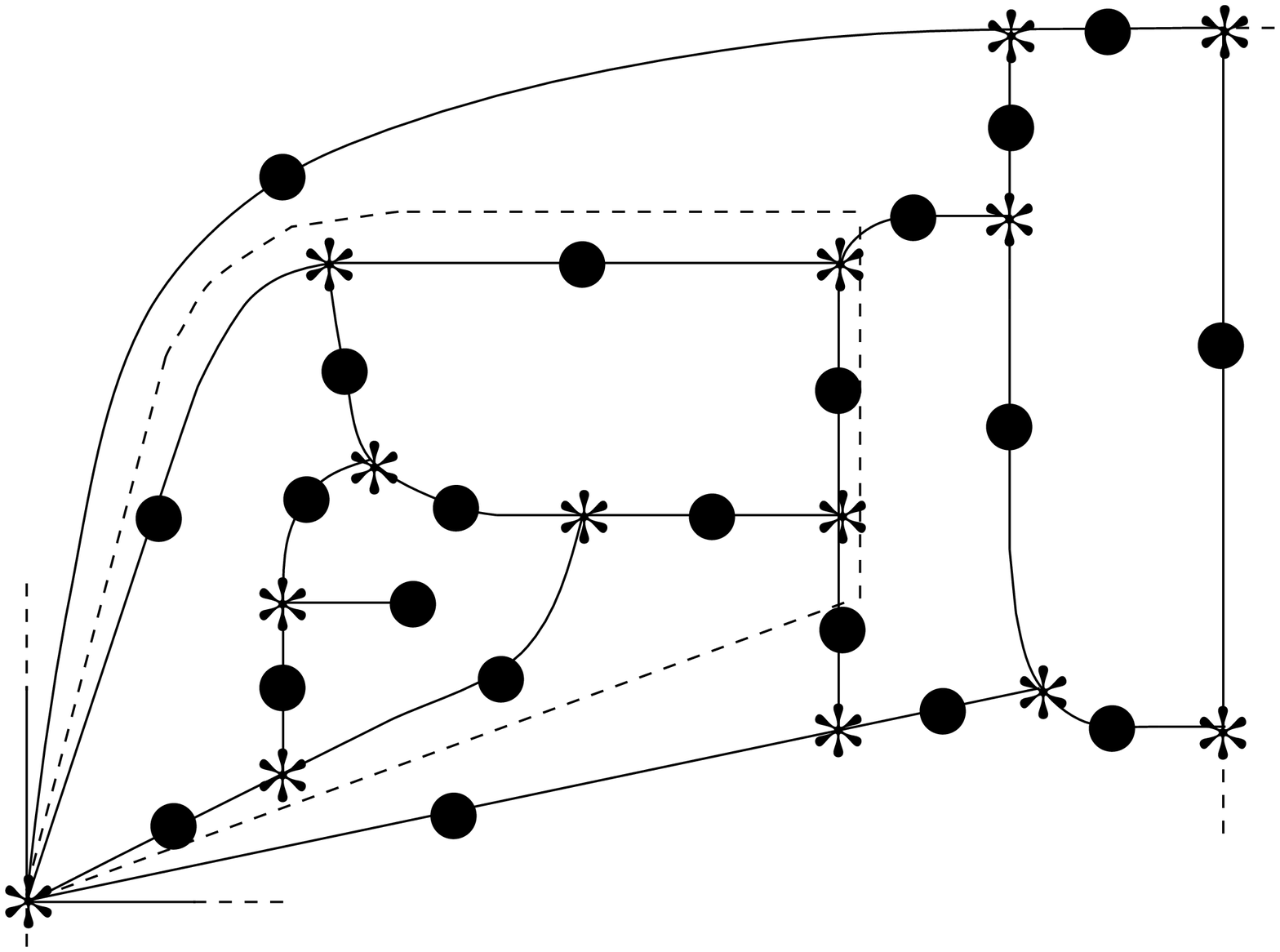}
\end{center}

We can easily see that the part which is enclosed by the dotted line 
in the above dessin repeatedly appears when $M$ is even 
and that dessins inductively constructed by this operation 
are compatible with the table. 

\vspace{2ex}
\underline{\textbf{(4) The case for $n \in \tfrac{1}{5}( \mathbb{N} + \tfrac{1}{2}$)}}

\begin{center}
  \begin{tabular}{|c|c|c|c|c|c|c|}                                                                                       \hline
     {}      & 0 & 1 & $\lambda$ & $\infty$          & $a_{1}, \cdots a_{12n-2}$                  & deg     \\ \hline 
     0       & 1 & 1 & 1         & 0                 & $( 15n-\tfrac{3}{2})$ pts with $mult. = 2$ & $30n$   \\ \hline
     1       & 0 & 0 & 0         & 0                 & $10n$ pts with $mult. = 3$                 & $30n$   \\ \hline
     $\infty$& 0 & 0 & 0         & $5n+\tfrac{5}{2}$ & $( 5n-\tfrac{1}{2})$ pts with $mult. = 5$                  & $30n$   \\ \hline
  \end{tabular}
\end{center}

If the dessin has $N$ loops with valency 5, 
\[
N = 5n - \tfrac{1}{2} \Leftrightarrow n = \tfrac{1}{5}(N + \tfrac{1}{2}) \in 
\tfrac{1}{5}(\mathbb{N} + \tfrac{1}{2}) 
\] 
and then the table becomes the following. 

\begin{center}
  \begin{tabular}{|c|c|c|c|c|c|c|}                                                                                       \hline
     {}      & 0 & 1 & $\lambda$ & $\infty$ & $a_{1}, \cdots a_{12n-2}$       & deg      \\ \hline 
     0       & 1 & 1 & 1         & 0        & $ 3N $ pts with $mult. = 2$     & $6N + 3$ \\ \hline
     1       & 0 & 0 & 0         & 0        & $(2N + 1)$ pts with $mult. = 3$ & $6N + 3$ \\ \hline
     $\infty$& 0 & 0 & 0         & $N + 3$  & $ N $ pts with $mult. = 5$      & $6N + 3$ \\ \hline
  \end{tabular}
\end{center}
So it suffices to construct the dessins compatible with the table 
for all $N \in \mathbb{N}$. 
Now, we construct dessins. 
(For each $N \geq 1$, our dessin is one of those which are compatible with the table.) 

\vspace{1ex}
For $N=0$, 
\begin{center}
\includegraphics[width=10em,clip]{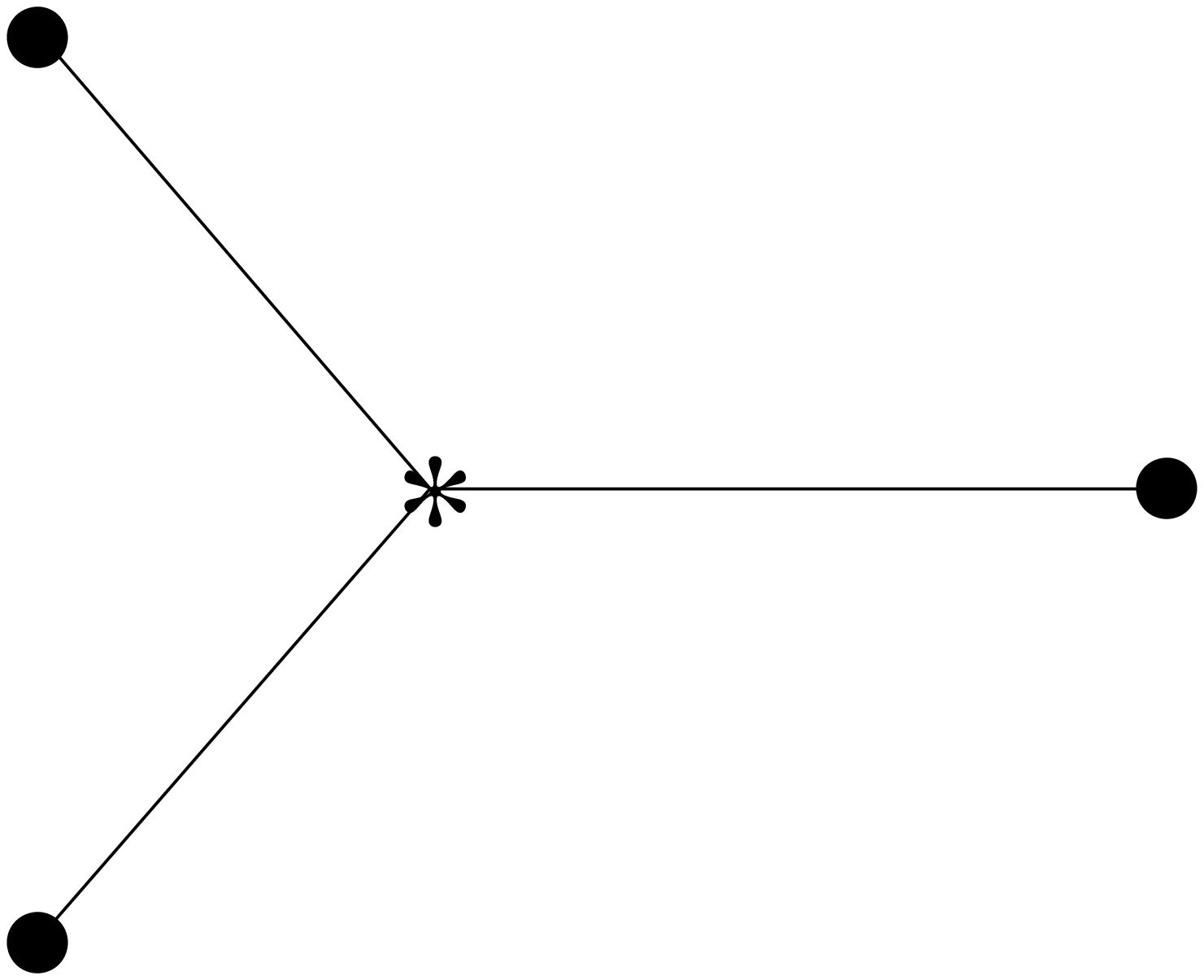}
\end{center}

For $N=1$, 
\begin{center}
\includegraphics[width=10em,clip]{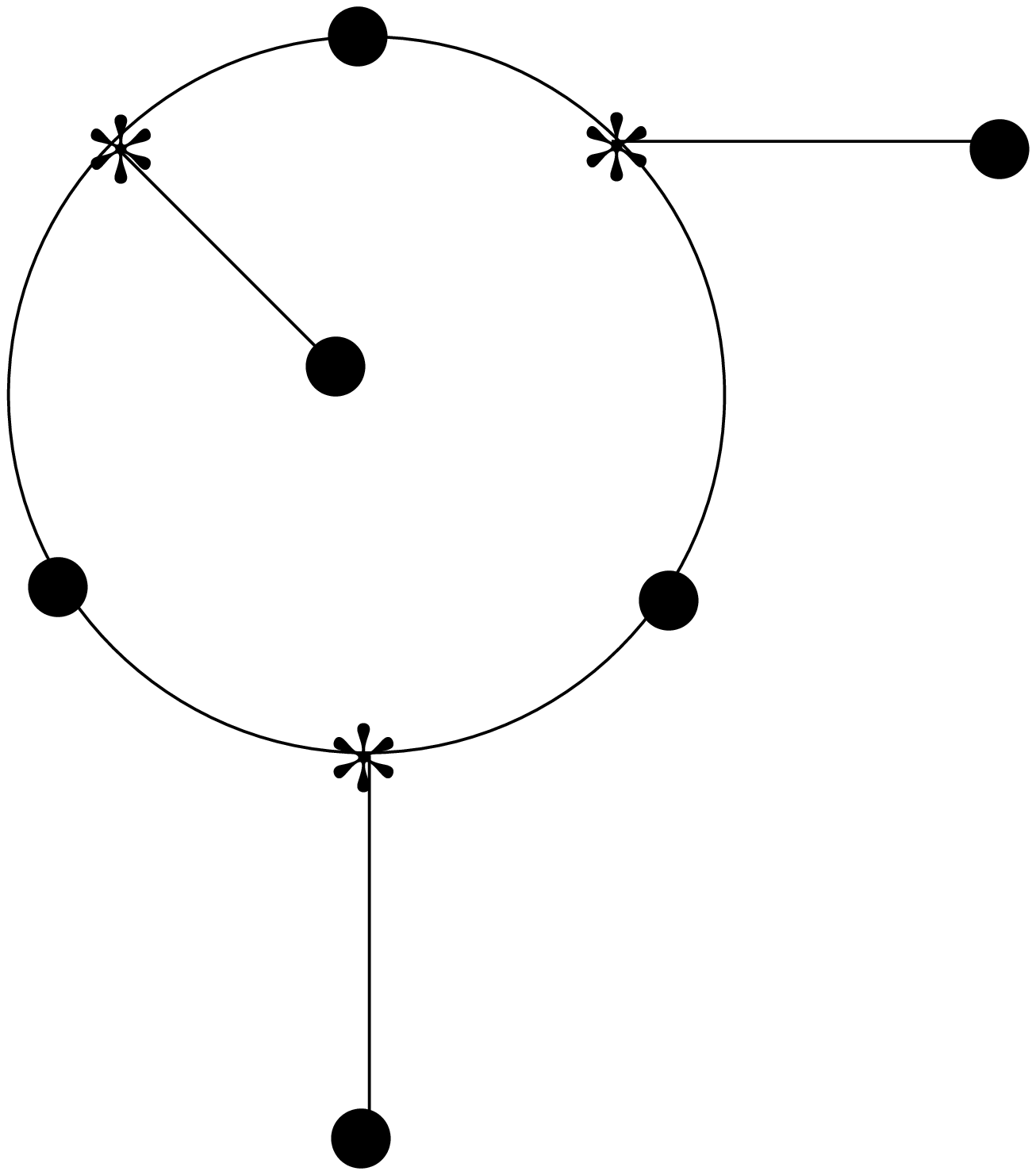}
\end{center}

For $N=2$, 
\begin{center}
\includegraphics[width=10em,clip]{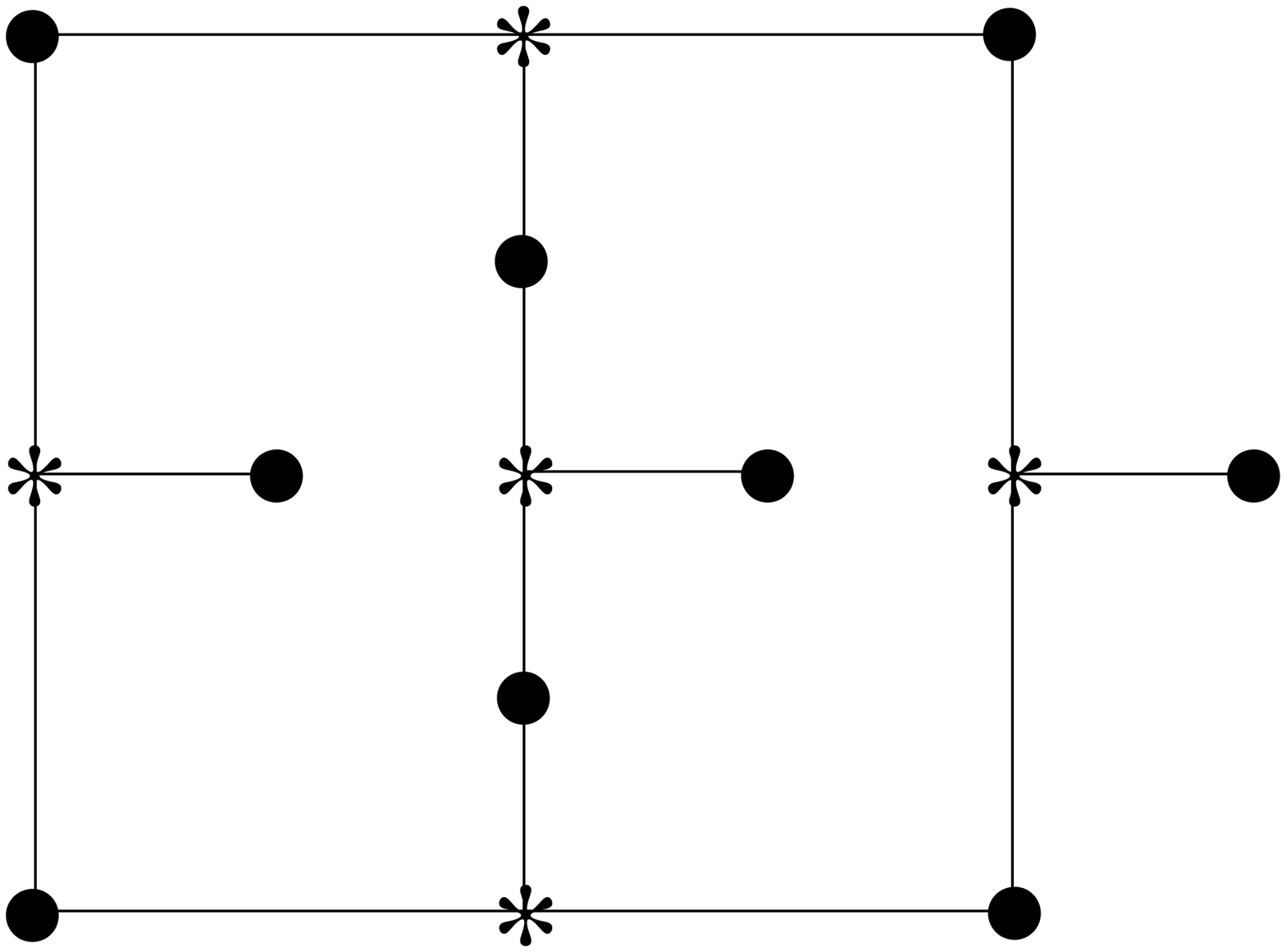}
\end{center}

For $N=3$, 
\begin{center}
\includegraphics[width=10em,clip]{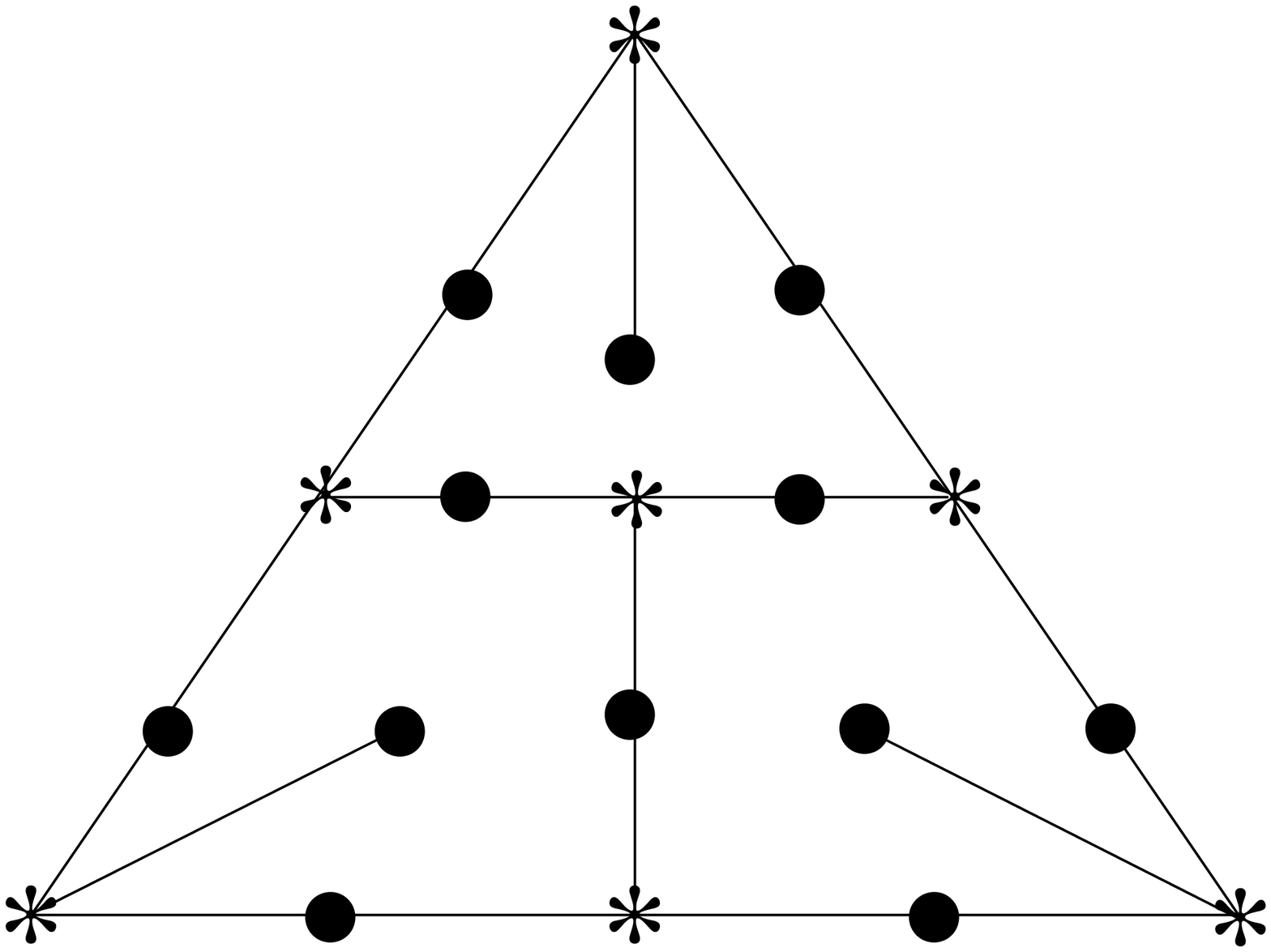}
\end{center}

For $N=4$, 
\begin{center}
\includegraphics[width=12em,clip]{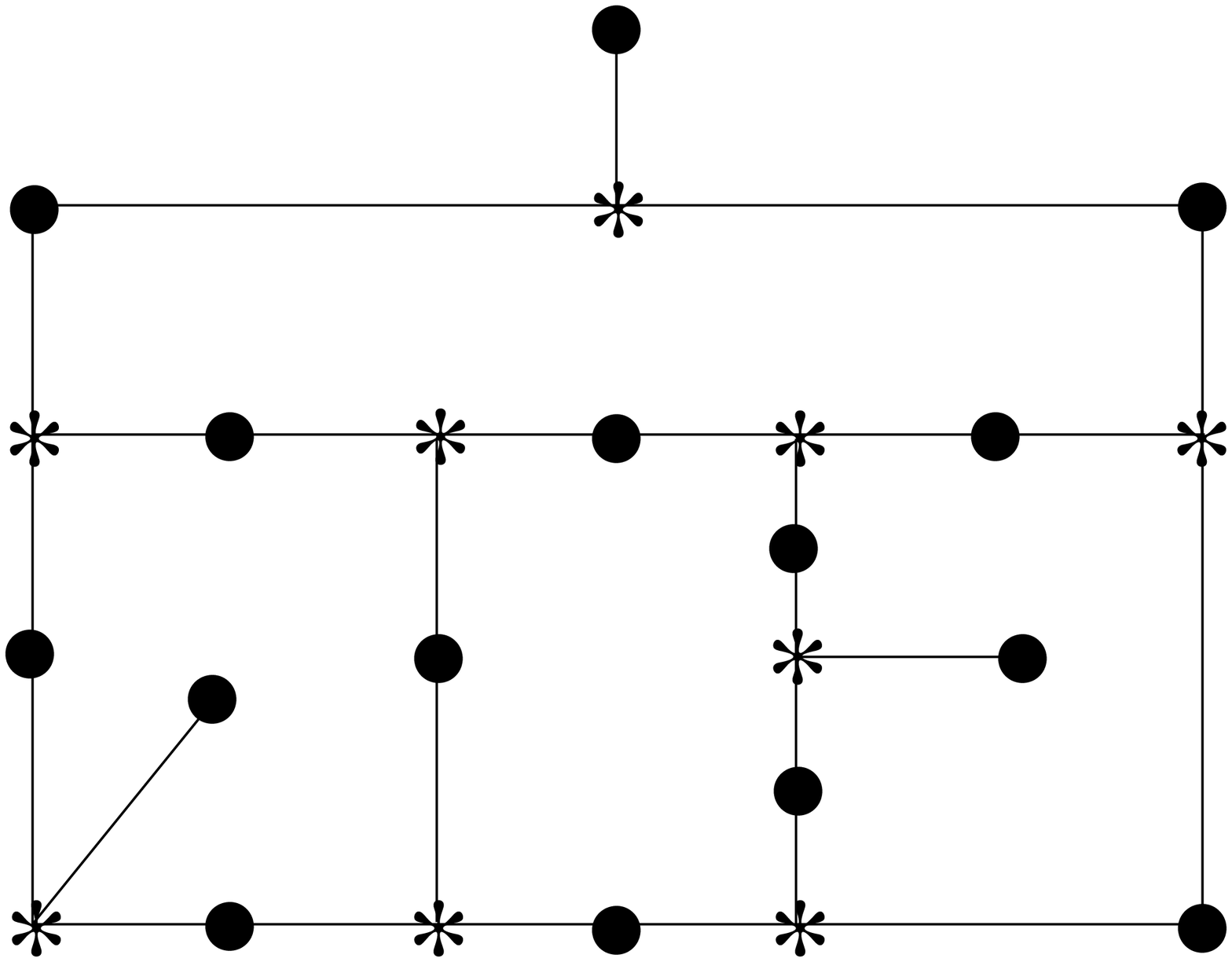}
\end{center}

For $N=5$, 
\begin{center}
\includegraphics[width=12em,clip]{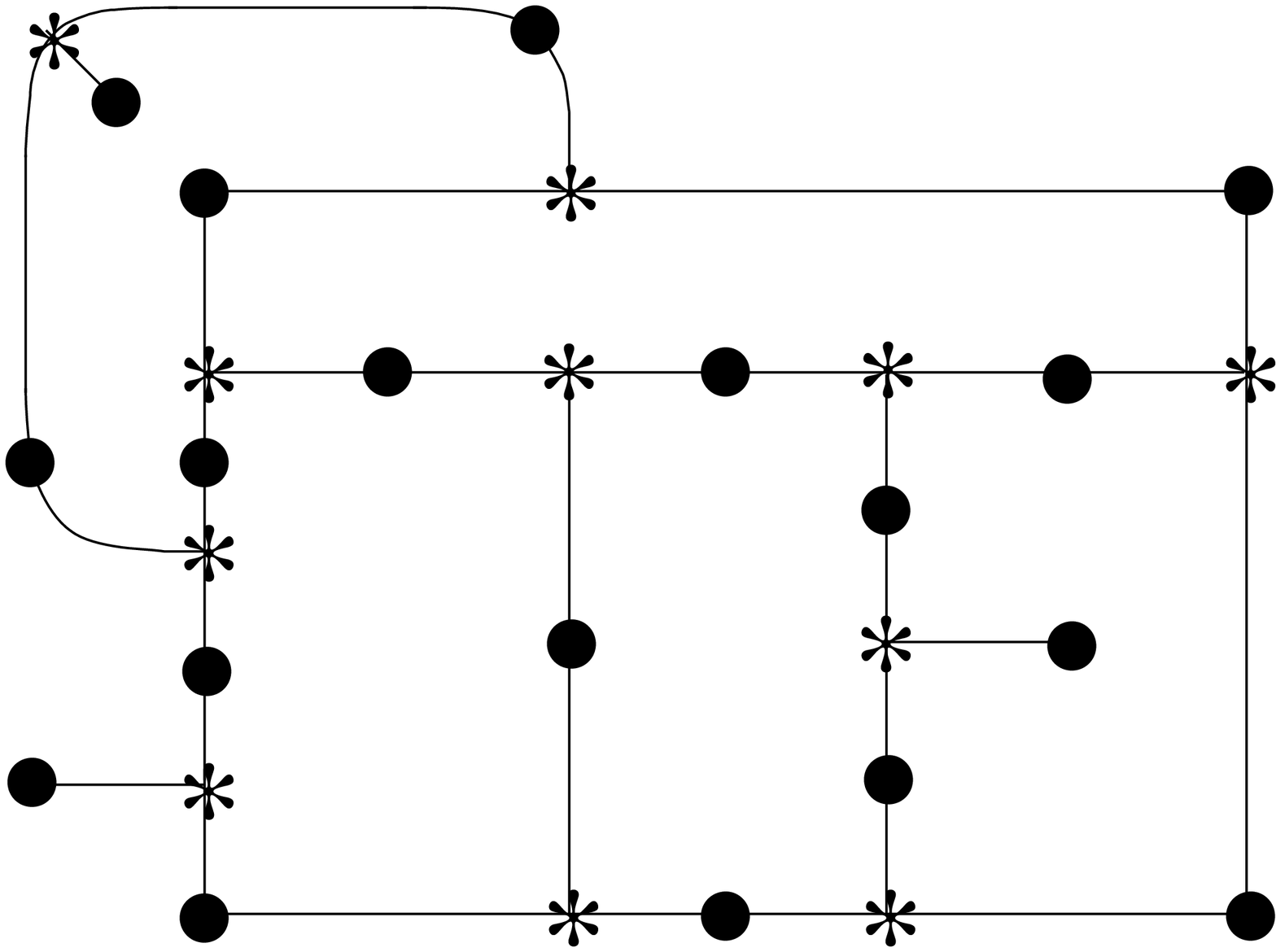}
\end{center}

For $N=6$, 
\begin{center}
\includegraphics[width=15em,clip]{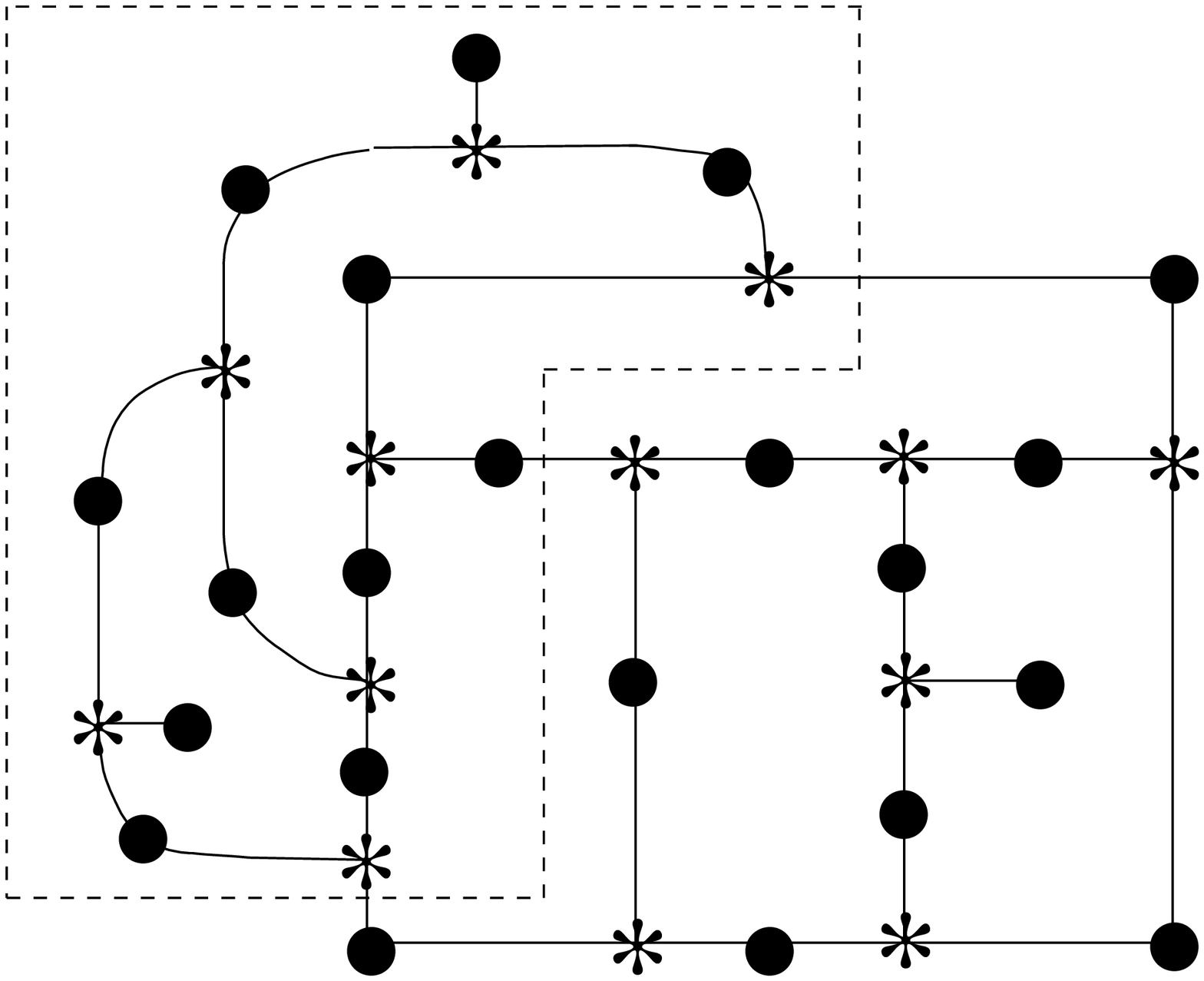}
\end{center}

\vspace{1ex}
For $N \geq 6$, we can construct dessins inductively according to 
the following operations. We operate the left side of the dessin. 
(We draw only the main part which is enclosed by the dotted line in the dessin of $N=6$.) 

If the dessin with $N=2k, (k \geq 3)$ loops is 
 
\begin{center}
\includegraphics[width=15em,clip]{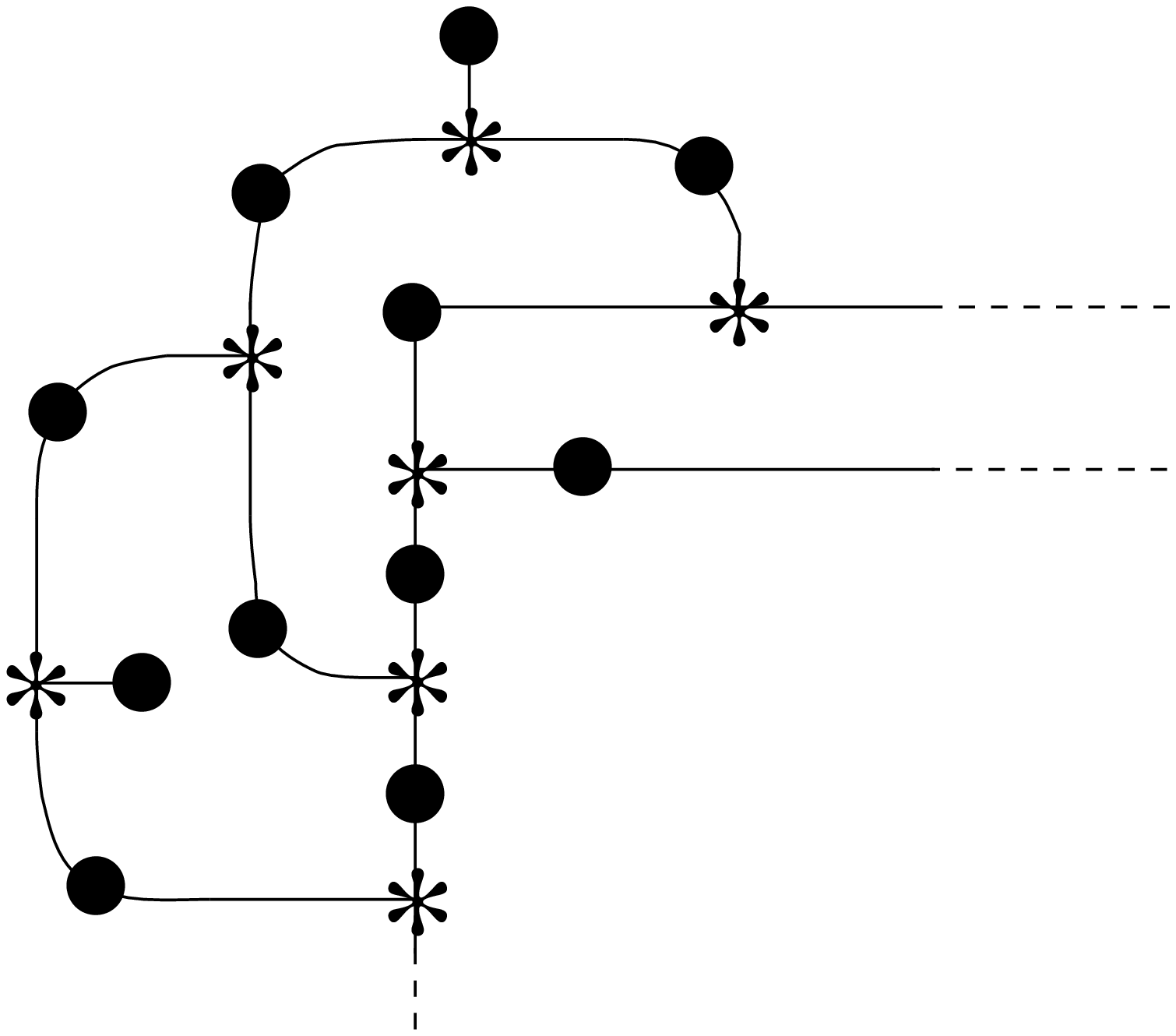}
\end{center}

\noindent
then the dessin with $N+1=2k+1$ loops is obtained by  

\begin{center}
\includegraphics[width=15em,clip]{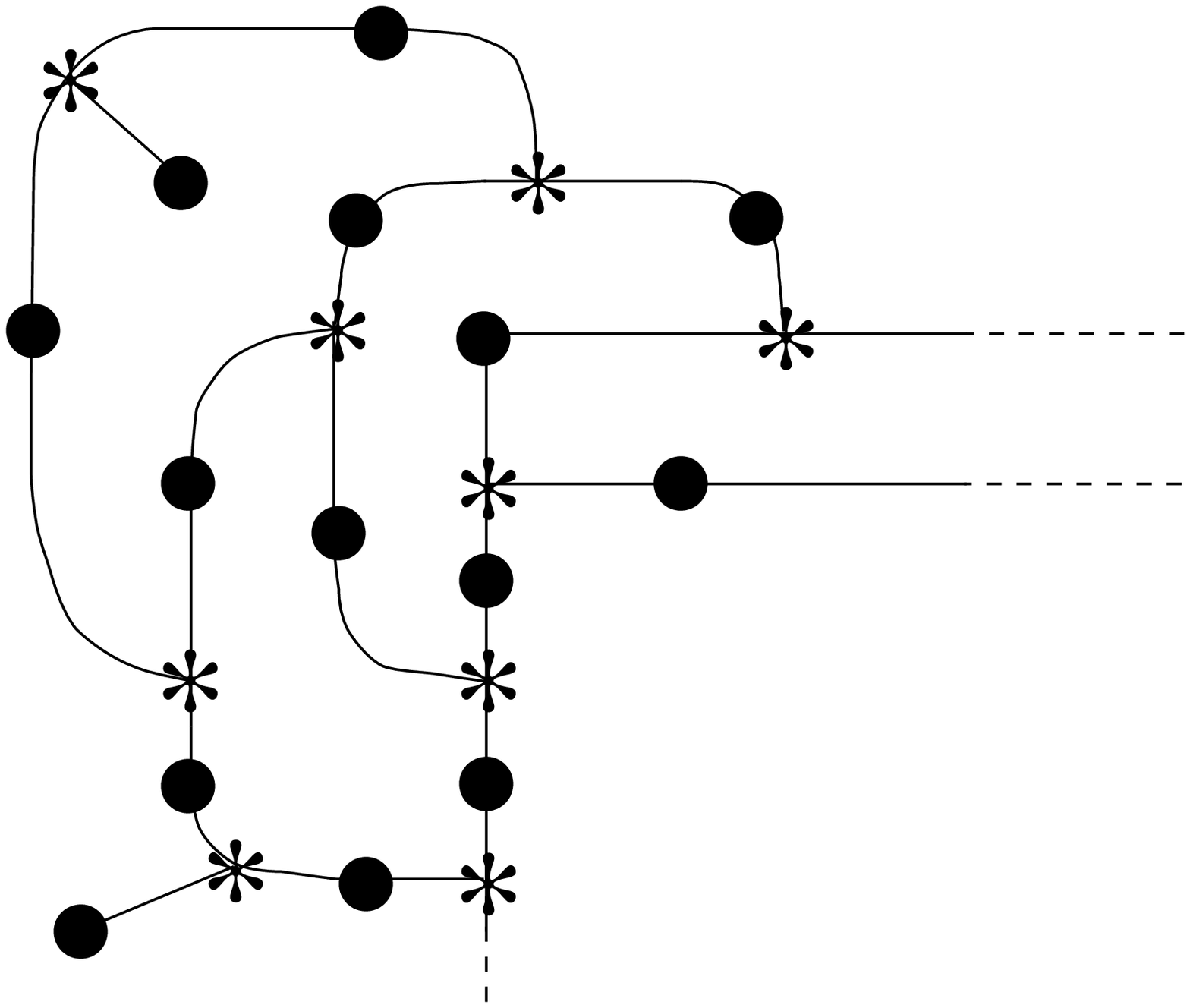}
\end{center}

\noindent
and the dessin with $N+2=2k+2$ loops is obtained by 
 
\begin{center}
\includegraphics[width=15em,clip]{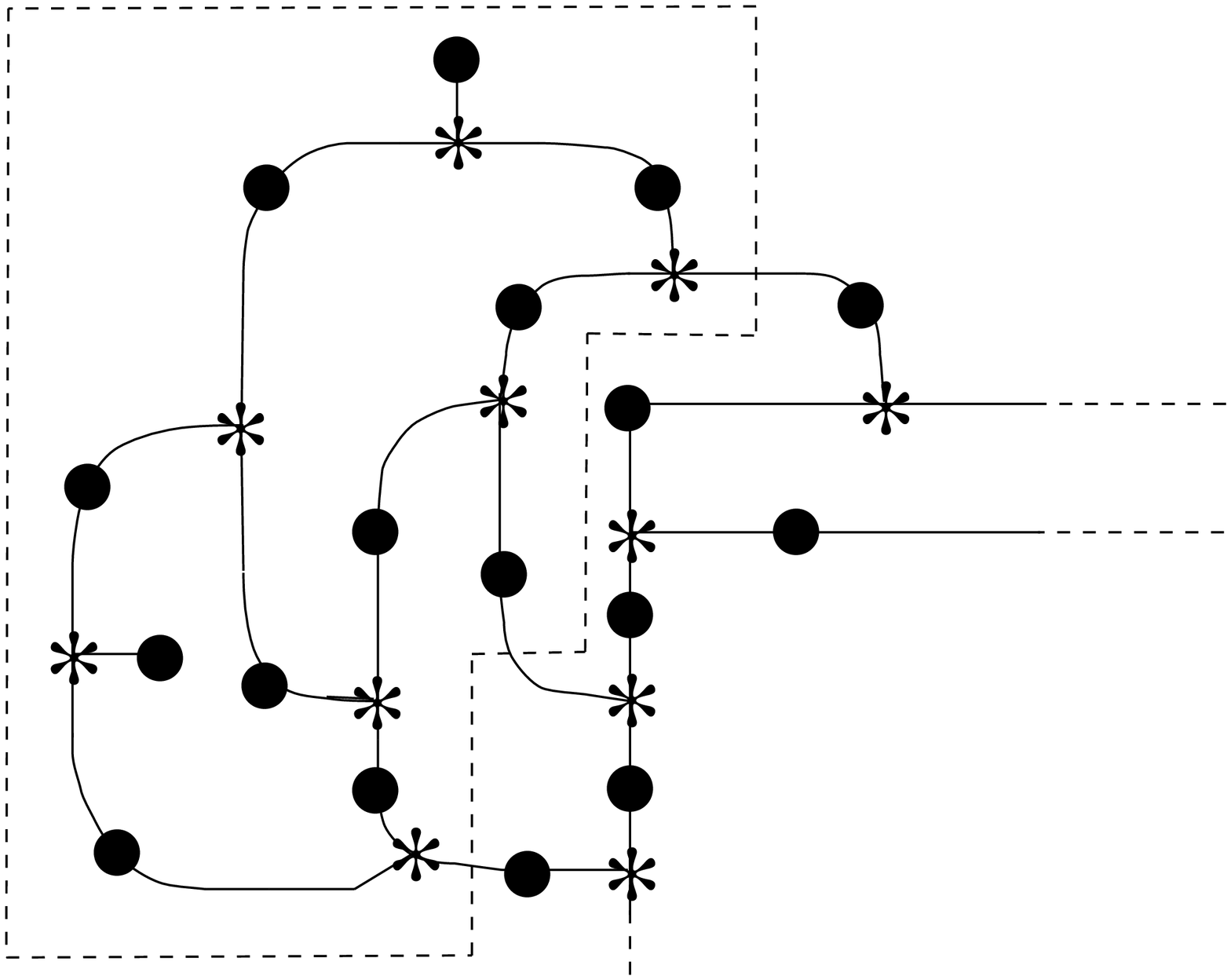}
\end{center}

We can easily see that the part which is enclosed by the dotted line 
in the above dessin repeatedly appears when $N$ is even 
and that dessins inductively constructed by this operation 
are compatible with the table. 

\begin{remark}\label{rem-exc2}
\rm{
The ``exceptional'' values of $n$ in icosahedral case are $n=-1/6$ and $n=-1/10$.
Similarly to the case of octahedral monodromy, we can show that these cases do not occur.
}
\end{remark}

\vspace{1ex}
\textbf{Acknowledgements.}
I would like to thank R\u{a}zvan Li\c{t}canu for showing me the new version of~\cite{L2}. 
I am deeply grateful to my adviser, Professor Fumiharu Kato for his advice and support.
Thanks are also due to Masao Aoki and Koki Itoh for their useful comments.

\noindent
Department of Mathematics, Faculty of Science, Kyoto University, Kyoto 606-8502, JAPAN.

\noindent
E-mail: \texttt{nakanisi@math.kyoto-u.ac.jp}

\end{document}